# Convergence rates of Forward–Douglas–Rachford splitting method

Cesare Molinari*, Jingwei Liang† and Jalal Fadili‡

**Abstract.** Over the past years, operator splitting methods have become ubiquitous for non-smooth optimization owing to their simplicity and efficiency. In this paper, we consider the Forward–Douglas–Rachford splitting method (FDR) [10, 40], and study both global and local convergence rates of this method. For the global rate, we establish an $o(1/k)$ convergence rate in terms of a Bregman divergence suitably designed for the objective function. Moreover, when specializing to the case of Forward–Backward splitting method, we show that convergence rate of the objective function of the method is actually $o(1/k)$ for a large choice of the descent step-size. Then locally, based on the assumption that the non-smooth part of the optimization problem is partly smooth, we establish local linear convergence of the method. More precisely, we show that the sequence generated by FDR method first (i) identifies a smooth manifold in a finite number of iteration, and then (ii) enters a local linear convergence regime, which is for instance characterized in terms of the structure of the underlying active smooth manifold. To exemplify the usefulness of the obtained result, we consider several concrete numerical experiments arising from applicative fields including, for instance, signal/image processing, inverse problems and machine learning.

**Key words.** Forward–Douglas–Rachford splitting, Forward–Backward splitting, Bregman distance, Partial smoothness, Finite identification, Local linear convergence.

**AMS subject classifications.** 49J52, 65K05, 65K10.

# 1 Introduction

## 1.1 Problem statement

Splitting methods are iterative algorithms to solve inclusion problems and optimization problems by decoupling the original problem into subproblems that are fast and easy to solve. These schemes evaluate the individual operators, their resolvents, the linear operators, all separately at various points in the course of iteration, but never the resolvents of sums nor of composition by a linear operator. Since the first operator splitting method developed in the 70's for solving structured monotone inclusion problems, the class of splitting methods has been regularly enriched with increasingly sophisticated algorithms as the structure of problems to handle becomes more complex. We refer the reader to [4] and references therein for a through account of operator splitting methods.

Nowadays, operator splitting methods have become ubiquitous for solving problems arising from fields through science and engineering, such as inverse problems, signal/image processing and machine learning, among others. Numerous applications encountered in these fields may end up with solving certain types of optimization problems. In this paper, we consider the following one:

$$\min_{x\in\mathcal{H}} \big\{F(x)+R(x):\quad x\in V\big\}, \tag{1.1}$$

*Universidad Técnica Federico Santa María, Av. España 1680, Valparaíso, Chile. E-mail: cecio.molinari@gmail.com.
†DAMTP, University of Cambridge, UK. E-mail: jl993@cam.ac.uk.
‡Normandie Université, ENSICAEN, CNRS, GREYC, France. E-mail: Jalal.Fadili@ensicaen.fr.

where $\mathcal{H}$ is a real Hilbert space equipped with scalar product $\langle \cdot, \cdot \rangle$ and norm $\|\cdot\|$. We suppose the following assumptions :

(**A.1**) The function $R$ belongs to $\Gamma_0(\mathcal{H})$, with $\Gamma_0(\mathcal{H})$ being the set of proper convex and lower semi-continuous functions from $\mathcal{H}$ to the extended real line $]-\infty, +\infty]$.

(**A.2**) $F : \mathcal{H} \to \mathbb{R}$ is convex continuously differentiable with $\nabla F$ being $(1/\beta)$-Lipschitz continuous.

(**A.3**) The constraint set $V$ is a closed vector subspace of $\mathcal{H}$.

(**A.4**) The set of minimizers is not empty.

Typical examples of (1.1) can be found in the numerical experiment section.

### 1.2 Forward–Douglas–Rachford splitting method

When $V = \mathcal{H}$, problem (1.1) can be handled by the classical Forward–Backward (FB) splitting method [33], whose iteration, in its relaxed form, reads

$$x_{k+1} = (1-\lambda_k)x_k + \lambda_k \mathrm{prox}_{\gamma R}\big(x_k - \gamma \nabla F(x_k)\big), \tag{1.2}$$

where $\gamma \in ]0, 2\beta[$ is the step-size and $\lambda_k \in ]0, \frac{4\beta-\gamma}{2\beta}[$ is the relaxation parameter. The term $\mathrm{prox}_{\gamma R}$ is called the proximity operator of $\gamma R$ and is defined by

$$\mathrm{prox}_{\gamma R}(x) \stackrel{\text{def}}{=} \mathrm{argmin}_{u\in\mathcal{H}}\, \gamma R(u) + \frac{1}{2}\|u - x\|^2. \tag{1.3}$$

When $V$ is merely a subspace of $\mathcal{H}$, in principle we still can apply FB splitting method to solve (1.1). However, even if $\mathrm{prox}_{\gamma R}$ is very easy to compute, the proximity operator of $R + \iota_V$ in general may be rather difficult to calculate. Therefore, new splitting algorithms are needed, and one possible choice is the Forward–Douglas–Rachford splitting method [10] which will be presented shortly in the next. Let us first define $\mathrm{P}_V$ as the orthogonal projector onto the subspace $V$, and the function $G \stackrel{\text{def}}{=} F \circ \mathrm{P}_V$. Then (1.1) is, obviously, equivalent to

$$\min_{x\in\mathcal{H}} \big\{ \Phi_V(x) \stackrel{\text{def}}{=} G(x) + R(x) + \iota_V(x) \big\}. \tag{1.4}$$

**Remark 1.1.** From the assumption on $F$, we have that also $G$ is convex and continuously differentiable with $\nabla G = \mathrm{P}_V \circ \nabla F \circ \mathrm{P}_V$ being $(1/\beta_V)$-Lipschitz continuous (notice that $\beta_V \geq \beta$). The observation of using $G$ instead of $F$ to achieve a better Lipschitz condition was first considered in [17].

The iteration of FDR method for solving (1.4) reads

$$\begin{cases} u_{k+1} = \mathrm{prox}_{\gamma R}\big(2x_k - z_k - \gamma\nabla G(x_k)\big), \\ z_{k+1} = z_k + \lambda_k(u_{k+1} - x_k), \\ x_{k+1} = \mathrm{P}_V(z_{k+1}), \end{cases} \tag{1.5}$$

where $\gamma$ is the step-size and $\lambda_k$ is the relaxation parameter. Recall that, under the additional assumptions that $\gamma \in ]0, 2\beta_V[$, $\lambda_k \in ]0, \frac{4\beta_V-\gamma}{2\beta_V}[$ and $\sum_{k\in\mathbb{N}} \lambda_k\big(\frac{4\beta_V-\gamma}{2\beta_V} - \lambda_k\big) = +\infty$, the sequences $\{u_k\}_{k\in\mathbb{N}}$ and $\{x_k\}_{k\in\mathbb{N}}$ converge to a solution; see [10, Theorem 4.2].

In this paper, we consider a non-stationary version of (1.5), namely when the step-size $\gamma$ may change along the iterations. The method is described below in Algorithm 1.

**Algorithm 1:** Non-stationary Forward–Douglas–Rachford splitting

**Initial**: $k = 0, z_0 \in \mathcal{H}, x_0 = \mathrm{P}_V(z_0)$.

**repeat**

Let $\gamma_k \in ]0, 2\beta_V[$ and $\lambda_k \in ]0, \frac{4\beta_V - \gamma}{2\beta_V}[$:

$$\left\lfloor \begin{aligned} u_{k+1} &= \mathrm{prox}_{\gamma_k R}\big(2x_k - z_k - \gamma_k \nabla G(x_k)\big), \\ z_{k+1} &= z_k + \lambda_k(u_{k+1} - x_k), \\ x_{k+1} &= \mathrm{P}_V(z_{k+1}). \end{aligned} \right. \tag{1.6}$$

$k = k + 1$;

**until** *convergence*;

**Remark 1.2.** For global convergence, one can also consider an inexact version of (1.6) by incorporating additive errors in the computation of $u_k$ and $x_k$, though we do not elaborate more on this for the sake of local convergence analysis. One can consult [29] for more details on this aspect.

In the next, we suppose the following main assumption on the algorithm parameters:

(**A.5**) The sequence of the step-sizes $\{\gamma_k\}_{k\in\mathbb{N}}$ and the one of the relaxation parameters $\{\lambda_k\}_{k\in\mathbb{N}}$ verify:

- $0 < \underline{\gamma} \leq \gamma_k \leq \bar{\gamma} < 2\beta_V$ and $\gamma_k \to \gamma$ for some $\gamma \in [\underline{\gamma}, \bar{\gamma}]$;
- $\lambda_k \in ]0, \frac{4\beta_V - \gamma_k}{2\beta_V}[$;
- $\sum_{k\in\mathbb{N}} \lambda_k(\frac{4\beta_V - \gamma_k}{2\beta_V} - \lambda_k) = +\infty$;
- $\sum_{k\in\mathbb{N}} \lambda_k |\gamma_k - \gamma| < +\infty$.

Notice that, for the stationary case (*i.e.* for $\gamma_k$ constant), assumption (**A.5**) is equivalent to the conditions required in [10, Theorem 4.2] for the convergence of iteration (1.5). Moreover, to satisfy (**A.5**) in absence of relaxation (*i.e.* when the relaxation parameter is fixed to $\lambda_k \equiv 1$), the sequence of the step-sizes has just to verify $\gamma_k \in ]\underline{\gamma}, \bar{\gamma}[$ with $\sum_{k\in\mathbb{N}} |\gamma_k - \gamma| < +\infty$. On the other hand, in general, the summability assumption of $\{\lambda_k|\gamma_k - \gamma|\}_{k\in\mathbb{N}}$ in (**A.5**) is weaker than imposing it without $\lambda_k$. Indeed, following the discussion in [13, Remark 5.7], take $q \in ]0, 1]$, let $\theta = \frac{4\beta_V - \bar{\gamma}}{4\beta_V} > \frac{1}{2}$ and

$$\lambda_k = \theta - \sqrt{\theta - 1/(2k)} \ \text{ and } \ |\gamma_k - \gamma| = \frac{\theta + \sqrt{\theta - 1/(2k)}}{k^q}.$$

Then it can be verified that

$$\textstyle\sum_{k\in\mathbb{N}} |\gamma_k - \gamma| = +\infty, \ \sum_{k\in\mathbb{N}} \lambda_k |\gamma_k - \gamma| = \frac{1}{2k^{1+q}} < +\infty$$

and

$$\textstyle\sum_{k\in\mathbb{N}} \lambda_k \big(\frac{4\beta_V - \gamma_k}{2\beta_V} - \lambda_k\big) \geq \sum_{k\in\mathbb{N}} \lambda_k (2\theta - \lambda_k) = \sum_{k\in\mathbb{N}} \frac{1}{2k} = +\infty.$$

As revealed by its name, FDR recovers the Douglas–Rachford splitting method [18] when $F = 0$, and also the FB splitting method [33] when $V = \mathcal{H}$. FDR method is also closely related to other two operator splitting methods: generalized Forward–Backward splitting (GFB) [40] and Three-operators splitting (TOS) [17], which are discussed below.

**Generalized Forward–Backward splitting** Let $m > 0$ be a positive integer. Now for problem (1.1), let be $V = \mathcal{H}$ and suppose we have $m$ non-smooth functionals. The problem then becomes

$$\min_{x\in\mathcal{H}} \big\{\Phi_m(x) \stackrel{\text{def}}{=} F(x) + \textstyle\sum_{i=1}^m R_i(x)\big\}, \tag{1.7}$$

where for each $i = 1, ..., m$, we have $R_i \in \Gamma_0(\mathcal{H})$.

Similar to the situation of FDR algorithm, even if the proximity operator of each $R_i$ can be solved easily, the proximity of the sum of them can be intractable. In [40], the authors propose the GFB algorithm, which achieves the full splitting of the evaluation of the proximity operator of each $R_i$. Let $(\omega_i)_i \in ]0,1[^m$ such that $\sum_{i=1}^m \omega_i = 1$, choose $\gamma \in ]0, 2\beta[$ and $\lambda_k \in ]0, \frac{4\beta-\gamma}{2\beta}[$:

$$\left\{\begin{aligned} &\text{from } i = 1 \text{ to } m: \\ &\quad\left\lfloor\begin{aligned} u_{i,k+1} &= \operatorname{prox}_{\frac{\gamma}{\omega_i}R_i}\big(2x_k - z_{i,k} - \gamma\nabla F(x_k)\big) \\ z_{i,k+1} &= z_{i,k} + \lambda_k(u_{i,k+1} - x_k) \end{aligned}\right. \\ &x_{k+1} = \textstyle\sum_{i=1}^m \omega_i z_{i,k+1}. \end{aligned}\right. \tag{1.8}$$

We refer to [40] for more details of the GFB algorithm. Now define the product space $\boldsymbol{\mathcal{R}} \stackrel{\text{def}}{=} \mathcal{H} \times ... \times \mathcal{H}$, equipped with proper inner product and norm, the subspace $\boldsymbol{\mathcal{S}} \stackrel{\text{def}}{=} \{\boldsymbol{x} = (x_i)_{i=1,...,m} \in \boldsymbol{\mathcal{R}} : x_1 = ... = x_m\} \subset \boldsymbol{\mathcal{R}}$ and let the weights be $\omega_i = \frac{1}{m}, i = 1, ..., m$. Then it can be shown that GFB algorithm is equivalent to applying FDR to the following problem:

$$\min_{\boldsymbol{x}\in\boldsymbol{\mathcal{R}}} F\Big(\frac{1}{m}\textstyle\sum_{i=1}^m x_i\Big) + \sum_{i=1}^m R_i(x_i) + \iota_{\boldsymbol{\mathcal{S}}}(\boldsymbol{x}).$$

We refer to [10] for more connections between FDR and GFB.

**Three-operator splitting** When in problem (1.7) we choose specifically $m = 2$, it becomes

$$\min_{x\in\mathcal{H}} F(x) + R_1(x) + R_2(x). \tag{1.9}$$

Notice that problem (1.9) can be handled by GFB as it is only a special case of (1.7). In [17] the author proposed a splitting scheme which resembles FDR yet different: given $\gamma \in ]0, 2\beta[$ and $\lambda_k \in ]0, \frac{4\beta-\gamma}{2\beta}[$, the iteration of Three-operator splitting (TOS) reads as follows:

$$\left\{\begin{aligned} u_{k+1} &= \operatorname{prox}_{\gamma R_1}\big(2x_k - z_k - \gamma\nabla F(x_k)\big) \\ z_{k+1} &= z_k + \lambda_k(u_{k+1} - x_k), \\ x_{k+1} &= \operatorname{prox}_{\gamma R_2}(z_{k+1}). \end{aligned}\right. \tag{1.10}$$

It can be observed that the projection operator $\mathrm{P}_V$ of FDR is replaced by the proximity operator $\operatorname{prox}_{\gamma R_2}$. Though the difference is only for the update of $x_{k+1}$, as we will see in Section 4.4, the fixed-point operator of the two methods are quite different.

## 1.3 Previous works

Over the past years, convergence rate analysis of some operator splitting methods have been extensively studied, typically in terms of the objective values or the sequence residuals. Typical examples are the FB splitting method and its accelerated version [37, 38, 6, 11, 2]. There are however many splitting algorithms where studying convergence rate of the objective value is not easy or impossible, unless additional assumptions on the involved functions are assumed. This is typically the case for splitting algorithms that are not of descent-type. Thus designing an appropriate criterion (e.g. deriving from some Lyapunov function) is challenging (as we will see for the FDR algorithm in Section 3).

The convergence rate of the objective value for FDR algorithm has been studied in [16]. There, the author presented some ergodic and pointwise convergence rates on the objective value under different (more or less stringent) assumptions imposed on the function $R$ in the objective (1.1). Without any further assumptions other than (**A.1**)-(**A.5**), the author proved a pointwise $o(1/\sqrt{k})$ convergence rate on the criterion $F(x_k) + R(u_k) - \Phi_V(x^\star)$ in absolute value (see [16, Theorem 3.5]). However, beside the fact this rate suggests that FDR seems quite pessimistic (it suggests that FDR is as slow as subgradient descent), there is no non-negativity guarantee for such a criterion and the obtained rate is thus of a quite limited interest. Improving this rate on the objective value requires quite strong assumptions on $R$.

As far as local linear convergence of the sequence in absence of strong convexity is concerned, it has received an increasing attention in the past few years in the context of first-order proximal splitting methods. The key idea here is to exploit the geometry of the underlying objective around its minimizers. This has been done for instance in [28, 30, 31, 32] for the FB scheme, Douglas–Rachford splitting/ADMM and Primal–Dual splitting, under the umbrella of partial smoothness. The error bound property[1] , as highlighted in the seminal work of [34, 35], were used by several authors to study linear convergence of first-order descent-type algorithms, and in particular FB splitting, see e.g. [8, 45, 27, 19]. However, to the best of our knowledge, we are not aware of local linear convergence results for the FDR algorithm.

## 1.4 Main contributions

The main contributions of this paper can be summarized as follows.

### 1.4.1 Global convergence

**Convergence of the non-stationary FDR** In Section 3, we first prove the convergence of the newly proposed non-stationary FDR scheme (1.6). This is achieved by capturing non-stationarity as an error term. The proof exploits a general result on inexact and non-stationary Krasnosel'skiĭ-Mann fixed-point iteration developed in [29].

**Convergence rate of a Bregman divergence** We design a Bregman divergence $D(u_k)$ (see (3.4)) as a meaningful convergence criterion. Under the sole assumptions (**A.1**)-(**A.5**), we show a pointwise $o(1/(k+1))$ and ergodic $O(1/(k+1))$ convergence rates of this criterion (Theorem 3.6). More precisely, for $k \in \mathbb{N}$, we show that

$$\min_{0\le i\le k} D(u_i) = o(1/(k+1)) \ \text{ and } \ D(\bar{u}_k) = O(1/(k+1)) \quad \text{where} \quad \bar{u}_k = \frac{1}{k+1}\sum_{i=0}^{k} u_i.$$

When specializing the result to the unrelaxed FB splitting method, we show that the convergence rate is actually $\Phi(x_k) - \Phi(x^\star) = o(1/(k+1))$ even for the larger interval $]0, 2\beta[$ on the step-size $\gamma_k$.

[1] For the interplay between the error bound property, the Kurdyka-Łojasiewicz property, and the quadratic growth property, see [8, 19].

### 1.4.2 Local convergence

We then turn to analyzing the local convergence behaviour where now $\mathcal{H}$ is assumed finite-dimensional, i.e. $\mathcal{H} = \mathbb{R}^n$. Let $x^\star \in \mathrm{Argmin}(\Phi_V)$ be a global minimizer of (1.4).

**Finite time activity identification** Under a non-degeneracy condition (see (4.2)), and provided that $R$ is partly smooth relative to a $C^2$-smooth manifold $\mathcal{M}^R_{x^\star}$ near $x^\star$ (see Definition 2.10), we show that the FDR scheme has a finite time identification property. More precisely, the FDR generated sequence $\{u_k\}_{k\in\mathbb{N}}$, which converges to $x^\star$, will identify in finite time the manifold $\mathcal{M}^R_{x^\star}$ (Theorem 4.1).

**Local linear convergence** Capitalizing on this finite identification result, we first show that the globally non-linear iteration (1.6) locally linearizes along the identified manifold $\mathcal{M}^R_{x^\star}$. Then, we show that the convergence becomes locally linear. The rate of linear convergence is characterized precisely in terms of the Riemannian structure of $\mathcal{M}^R_{x^\star}$.

**Paper organization** The rest of the paper is organized as follows. In Section 2 we recall some classical material on convex analysis, operator theory that are essential to our exposition. We then introduce the notion of partial smoothness. The global convergence analysis is presented in Section 3, followed by finite identification and local convergence analysis in Section 4. Several numerical experiments are presented in Section 5. To keep the readability of the paper, very long proofs are collected in the appendix.

## 2 Preliminaries

Throughout the paper, $\mathcal{H}$ is a Hilbert space equipped with scalar product $\langle\cdot,\,\cdot\rangle$ and norm $\|\cdot\|$. $\mathrm{Id}$ denotes the identity operator on $\mathcal{H}$.

**Sets** For a nonempty convex set $C \subset \mathcal{H}$, denote by $\mathrm{aff}(C)$ its affine hull and by $\mathrm{par}(C)$ the smallest subspace parallel to $\mathrm{aff}(C)$. Denote $\iota_C$ the indicator function of $C$, $\mathcal{N}_C$ the associated normal cone operator and $\mathrm{P}_C$ the orthogonal projection operator on $C$.

**Functions** The sub-differential of a function $R \in \Gamma_0(\mathcal{H})$ is the set-valued operator defined by

$$\partial R : \mathcal{H} \rightrightarrows \mathcal{H},\; x \mapsto \big\{g \in \mathcal{H} | R(x') \geq R(x) + \langle g,\, x' - x\rangle, \forall x' \in \mathcal{H}\big\}. \tag{2.1}$$

**Lemma 2.1 (Descent lemma [7]).** *Suppose that $F : \mathcal{H} \to \mathbb{R}$ is convex continuously differentiable and $\nabla F$ is $(1/\beta)$-Lipschitz continuous. Then, given any $x, y \in \mathcal{H}$,*

$$F(x) \leq F(y) + \langle \nabla F(y),\, x - y\rangle + \frac{1}{2\beta}\|x - y\|^2.$$

**Definition 2.2 (Bregman divergence).** Given a function $R \in \Gamma_0(\mathcal{H})$ and two points $x, y$ in its effective domain, the Bregman divergence is defined by

$$\mathcal{D}^p_R(y, x) \overset{\text{def}}{=} R(y) - R(x) - \langle p,\, y - x\rangle,$$

where $p \in \partial R(x)$ is a sub-gradient of $R$.

Notice that the Bregman divergence is not a distance in the usual sense, since it is in general not symmetric[2]. However, it measures the distance of two points in the sense that $\mathcal{D}_R^p(x,x) = 0$ and $\mathcal{D}_R^p(y,x) \geq 0$ for any $x, y$ in the domain of $R$. Moreover, $\mathcal{D}_R^p(y,x) \geq \mathcal{D}_R^p(w,x)$ for all the points $w$ that are in the line segment between $x$ and $y$.

**Operators** Given a set-valued mapping $A : \mathcal{H} \rightrightarrows \mathcal{H}$, define its range $\mathrm{ran}(A) = \{y \in \mathcal{H} : \exists x \in \mathcal{H} \text{ s.t. } y \in A(x)\}$ and its graph as $\mathrm{gph}\,(A) \stackrel{\text{def}}{=} \{(x,u) \in \mathcal{H} \times \mathcal{H} :\ u \in A(x)\}$.

**Definition 2.3 (Monotone operator).** A set-valued mapping $A : \mathcal{H} \rightrightarrows \mathcal{H}$ is said to be monotone if,

$$\langle x_1 - x_2,\, v_1 - v_2 \rangle \geq 0,\ \forall (x_1, v_1) \in \mathrm{gph}\,(A) \ \text{ and } \ (x_2, v_2) \in \mathrm{gph}\,(A). \tag{2.2}$$

It is maximal monotone if $\mathrm{gph}\,(A)$ can not be contained in the graph of any other monotone operators.

For a maximal monotone operator $A$, $(\mathrm{Id} + A)^{-1}$ denotes its resolvent. It is known that for $R \in \Gamma_0(\mathcal{H})$, $\partial R$ is maximal monotone [41], and that $\mathrm{prox}_R = (\mathrm{Id} + \partial R)^{-1}$.

**Definition 2.4 (Cocoercive operator).** Let $\beta \in ]0, +\infty[$, $B : \mathcal{H} \to \mathcal{H}$, then $B$ is $\beta$-cocoercive if

$$\langle B(x_1) - B(x_2),\, x_1 - x_2 \rangle \geq \beta \|B(x_1) - B(x_2)\|^2,\ \forall x_1, x_2 \in \mathcal{H}. \tag{2.3}$$

If an operator is $\beta$-cocoercive, then it is $\beta^{-1}$-Lipschitz continuous.

**Definition 2.5 (Non-expansive operator).** An operator $\mathscr{F} : \mathcal{H} \to \mathcal{H}$ is non-expansive if

$$\|\mathscr{F}(x) - \mathscr{F}(y)\| \leq \|x - y\|, \quad \forall x, y \in \mathcal{H}.$$

That is, $\mathscr{F}$ is 1-Lipschitz continuous. For any $\alpha \in ]0,1[$, $\mathscr{F}$ is called $\alpha$-averaged if there exists a non-expansive operator $\mathscr{F}'$ such that $\mathscr{F} = \alpha\mathscr{F}' + (1-\alpha)\mathrm{Id}$.

In particular, when $\alpha = \frac{1}{2}$, $\mathscr{F}$ is called *firmly non-expansive*. Several properties of firmly non-expansive operators are collected in the following lemma.

**Lemma 2.6.** *Let $\mathscr{F} : \mathcal{H} \to \mathcal{H}$, the following statements are equivalent:*

(i) *$\mathscr{F}$ is firmly non-expansive;*

(ii) *$\mathrm{Id} - \mathscr{F}$ is firmly non-expansive;*

(iii) *$2\mathscr{F} - \mathrm{Id}$ is non-expansive;*

(iv) *$\mathscr{F}$ is the resolvent of a maximal monotone operator $A : \mathcal{H} \rightrightarrows \mathcal{H}$.*

**Proof.** (i)⇔(ii)⇔(iii) follows [4, Proposition 4.2, Corollary 4.29], and (i)⇔(iv) is [4, Corollary 23.8]. □

Some properties of functions with Lipschitz continuous gradient are shown in the next lemma.

**Lemma 2.7.** *Let $F : \mathcal{H} \to ]-\infty, +\infty[$ be a convex differentiable function, with $\frac{1}{\beta}$-Lipschitz continuous gradient, $\beta \in ]0, +\infty[$, then*

(i) *$\beta\nabla F$ is firmly non-expansive.*

(ii) *$\mathrm{Id} - \gamma\nabla F$ is $\frac{\gamma}{2\beta}$-averaged for $\gamma \in ]0, 2\beta[$.*

**Proof.** (i) See [3, Baillon–Haddad theorem]; (ii) See [4, Chapter 4, Proposition 4.33]; □

[2] It is symmetric if and only if $R$ is a nondegenerate convex quadratic form.

The class of $\alpha$-averaged operators is closed under relaxation, convex combination and composition [4]. The next lemma shows the composition of two averaged non-expansive operators..

**Lemma 2.8 ([39, Theorem 3]).** *Let $\mathscr{F}_1 : \mathcal{H} \to \mathcal{H}$ be $\alpha_1$-averaged and $\mathscr{F}_2 : \mathcal{H} \to \mathcal{H}$ be $\alpha_2$-average, then the composite operator $\mathscr{F}_1 \circ \mathscr{F}_2$ is $\alpha$-averaged with $\alpha = \frac{\alpha_1+\alpha_2-2\alpha_1\alpha_2}{1-\alpha_1\alpha_2} \in ]0, 1[$.*

The next lemma shows some properties of the set of fixed points for a non-expansive operator.

**Lemma 2.9.** *Let $\mathscr{F} : \mathcal{H} \to \mathcal{H}$ be non-expansive, then the set of fixed points $\mathrm{fix}(\mathscr{F}) = \{z \in \mathcal{H} | z = \mathscr{F}(z)\}$ is closed and convex.*

**Proof.** See [4, Chapter 4, Corollary 4.15]. □

**Partial smoothness** In this part, we assume that $\mathcal{H} = \mathbb{R}^n$ and we briefly introduce the concept of partial smoothness, on which relies our local convergence rate analysis. Partial smoothness was first introduced in [26]. This concept, as well as that of identifiable surfaces [44], captures the essential features of the geometry of non-smoothness when it can be localized along the so-called active/identifiable manifold.

Let $\mathcal{M}$ be a $C^2$-smooth embedded submanifold of $\mathbb{R}^n$ around a point $x$. To lighten the notation, henceforth we write $C^2$-manifold instead of $C^2$-smooth embedded submanifold of $\mathbb{R}^n$. The natural embedding of a submanifold $\mathcal{M}$ into $\mathbb{R}^n$ permits to define a Riemannian structure on $\mathcal{M}$, and we simply say that $\mathcal{M}$ is a Riemannian manifold. $\mathcal{T}_{\mathcal{M}}(x)$ denotes the tangent space to $\mathcal{M}$ at any point near $x$ in $\mathcal{M}$. More materials on manifolds are given in Section B.1.

Below we present the definition of partly smooth functions in $\Gamma_0(\mathbb{R}^n)$ setting.

**Definition 2.10 (Partly smooth function).** Let $R \in \Gamma_0(\mathbb{R}^n)$, and $x \in \mathbb{R}^n$ such that $\partial R(x) \neq \emptyset$. $R$ is then said to be *partly smooth* at $x$ relative to a set $\mathcal{M}$ containing $x$ if

(i) **Smoothness**: $\mathcal{M}$ is a $C^2$-manifold around $x$, $R$ restricted to $\mathcal{M}$ is $C^2$ around $x$;
(ii) **Sharpness**: The tangent space $\mathcal{T}_{\mathcal{M}}(x)$ coincides with $T_x \stackrel{\text{def}}{=} \mathrm{par}(\partial R(x))^{\perp}$;
(iii) **Continuity**: The set-valued mapping $\partial R$ is continuous at $x$ relative to $\mathcal{M}$.

The class of partly smooth functions at $x$ relative to $\mathcal{M}$ is denoted as $\mathrm{PSF}_x(\mathcal{M})$.

Under transversality assumptions, the set of partly smooth functions is closed under addition and pre-composition by a linear operator [26]. Moreover, absolutely permutation-invariant convex and partly smooth functions of the singular values of a real matrix, *i.e.* spectral functions, are convex and partly smooth spectral functions of the matrix [15]. Popular examples of partly smooth functions are summarized in Section 5 whose details can be found in [30].

The next lemma gives expressions of the Riemannian gradient and Hessian (see Section B.1 for definitions) of a partly smooth function.

**Lemma 2.11.** *If $R \in \mathrm{PSF}_x(\mathcal{M})$, then for any $x' \in \mathcal{M}$ near $x$,*

$$\nabla_{\mathcal{M}} R(x') = \mathrm{P}_{T_{x'}}(\partial R(x')).$$

*In turn, for all $h \in T_{x'}$,*

$$\nabla^2_{\mathcal{M}} R(x')h = \mathrm{P}_{T_{x'}} \nabla^2 \widetilde{R}(x')h + \mathfrak{W}_{x'}\big(h, \mathrm{P}_{T_{x'}^{\perp}} \nabla \widetilde{R}(x')\big),$$

*where $\widetilde{R}$ is any smooth extension (representative) of $R$ on $\mathcal{M}$, and $\mathfrak{W}_x(\cdot,\cdot) : T_x \times T_x^{\perp} \to T_x$ is the Weingarten map of $\mathcal{M}$ at $x$.*

**Proof.** See [30, Fact 3.3]. □

# 3 Global convergence

In this section, we deliver the global convergence behaviour of the non-stationary FDR algorithm (1.6) in a general real Hilbert space setting, including rate analysis.

## 3.1 Global convergence of the non-stationary FDR

Define the reflection operators of $\gamma R$ and $\iota_V$ respectively as $\mathcal{R}_{\gamma R} \stackrel{\text{def}}{=} 2\mathrm{prox}_{\gamma R} - \mathrm{Id}$ and $\mathcal{R}_V \stackrel{\text{def}}{=} 2\mathrm{P}_V - \mathrm{Id}$. Moreover, define the following operators:

$$\mathscr{F}_\gamma = \frac{1}{2}(\mathrm{Id} + \mathcal{R}_{\gamma R} \circ \mathcal{R}_V)(\mathrm{Id} - \gamma \nabla G) \quad \text{and} \quad \mathscr{F}_{\gamma,\lambda_k} = (1-\lambda_k)\mathrm{Id} + \lambda_k \mathscr{F}_\gamma. \tag{3.1}$$

The next lemma shows the property of the fixed-point operator $\mathscr{F}_{\mathrm{FDR}}$ of FDR algorithm.

**Lemma 3.1.** *For the FDR algorithm (1.6), let $\gamma \in ]0, 2\beta_V[$ and $\lambda_k \in ]0, \frac{4\beta_V - \gamma}{2\beta_V}[$. Then we have that $\mathscr{F}_\gamma$ is $\frac{2\beta_V}{4\beta_V-\gamma}$-averaged and $\mathscr{F}_{\gamma,\lambda_k}$ is $\frac{2\beta_V\lambda_k}{4\beta_V-\gamma}$-averaged.*

**Proof.** The averaged property of $\mathscr{F}_\gamma$ is a combination of [10, Proposition 4.1] and Lemma 2.8. For $\mathscr{F}_{\gamma,\lambda_k}$, it is sufficient to apply the definition of averaged operators. ☐

The (stationary) FDR iteration (1.5) can be written into a fixed-point iteration in terms of $z_k$ [10, Theorem 4.2], namely

$$z_{k+1} = \mathscr{F}_{\gamma,\lambda_k}(z_k). \tag{3.2}$$

As mentioned in the introduction, from [10, Theorem 4.2] we have the following convergence result. If $\mathrm{Argmin}(\Phi_V) \neq \emptyset$ and $\sum_{k\in\mathbb{N}} \lambda_k (\frac{4\beta_V-\gamma}{2\beta_V} - \lambda_k) = +\infty$, then $\{z_k\}_{k\in\mathbb{N}}$ converges weakly to some $z^\star \in \mathrm{fix}(\mathscr{F}_\gamma)$, and the sequence $\{x_k\}_{k\in\mathbb{N}}$ converges weakly to $x^\star \stackrel{\text{def}}{=} \mathrm{P}_V(z^\star)$, where $\mathrm{P}_V(z^\star) \in \mathrm{Argmin}(\Phi_V)$. On the other hand, the non-stationary FDR iteration (1.6) can be written as

$$\begin{aligned} z_{k+1} = \mathscr{F}_{\gamma_k,\lambda_k}(z_k) &= (1-\lambda_k)z_k + \lambda_k \mathscr{F}_{\gamma_k}(z_k) \\ &= \big((1-\lambda_k)z_k + \lambda_k \mathscr{F}_\gamma(z_k)\big) + \lambda_k\big(\mathscr{F}_{\gamma_k}(z_k) - \mathscr{F}_\gamma(z_k)\big). \end{aligned} \tag{3.3}$$

Then we have the following result about global convergence of the algorithm.

**Theorem 3.2 (Global convergence of non-stationary FDR).** *Consider the non-stationary FDR iteration (1.6). Suppose that Assumptions (A.1)-(A.5) hold. Then $\sum_{k\in\mathbb{N}} \|z_k - z_{k-1}\|^2 < +\infty$. Moreover, the sequence $\{z_k\}_{k\in\mathbb{N}}$ converges weakly to a point $z^\star \in \mathrm{fix}(\mathscr{F}_\gamma)$, and $\{x_k\}_{k\in\mathbb{N}}$ converges weakly to $x^\star \stackrel{\text{def}}{=} \mathrm{P}_V(z^\star) \in \mathrm{Argmin}(\Phi_V)$. If, in addition, either $\inf_{k\in\mathbb{N}} \lambda_k > 0$ or $\mathcal{H}$ is finite-dimensional, then $\{u_k\}_{k\in\mathbb{N}}$ also converges weakly to $x^\star$.*

The main idea of the proof of the theorem (see Section A) is to treat the non-stationarity as a perturbation error of the stationary iteration.

**Remark 3.3.**

(i) As emphasized in the introduction, Theorem 3.2 remains true if the iteration is carried out approximately, *i.e.* if $\mathscr{F}_{\gamma_k}(z_k)$ is computed approximately, provided that the errors are summable; See [29, Section 6] for more details.

(ii) With more assumptions on how fast $\{\gamma_k\}_{k\in\mathbb{N}}$ converges to $\gamma$, we can also derive the convergence rate of the sequence of residuals $\{\|z_k - z_{k-1}\|\}_{k\in\mathbb{N}}$. However, as we will study in Section 4 local linear convergence behaviour of $\{z_k\}_{k\in\mathbb{N}}$, we shall forgo the discussion here. Interested readers can consult [29] for more details about the rate of residuals.

### 3.2 Convergence rate of the Bregman divergence

As in Theorem 3.2, let $z^\star \in \mathrm{fix}(\mathscr{F}_\gamma)$ be a fixed point of (1.6) and $x^\star \stackrel{\text{def}}{=} \mathrm{P}_V(z^\star) \in \mathrm{Argmin}(\Phi_V)$. Then the following optimality condition holds

$$v^\star \in \partial R(x^\star) + \nabla G(x^\star),$$

where $v^\star \in V^\perp \stackrel{\text{def}}{=} \mathcal{N}_V(x^\star)$ is a normal vector. Now denote $\Phi \stackrel{\text{def}}{=} R + G$. For $y \in \mathbb{R}^n$, define the following Bregman divergence

$$\begin{aligned} \mathcal{D}_\Phi^{v^\star}(y) \stackrel{\text{def}}{=} \mathcal{D}_\Phi^{v^\star}(y, x^\star) &= \Phi(y) - \Phi(x^\star) - \langle v^\star,\, y - x^\star\rangle \\ &= \Phi(y) - \Phi(x^\star) - \langle v^\star,\, y^{V^\perp}\rangle, \end{aligned} \tag{3.4}$$

where $y^{V^\perp} \stackrel{\text{def}}{=} \mathrm{P}_{V^\perp}(y)$ is the projection of $y$ onto $V^\perp$, and we use the fact that $\langle v^\star,\, x^\star\rangle = 0$.

The motivation of choosing the above function to quantify the convergence rate of FDR algorithm is due to the fact that it measures both the discrepancy of the objective to the optimal value and violation of the constraint on $V$.

Lemma 3.4 hereafter will provide us with a key estimate on $\mathcal{D}_\Phi^{v^\star}(u_k)$ which in turn will be instrumental to derive the convergence rate of $\{\mathcal{D}_\Phi^{v^\star}(u_k)\}_{k\in\mathbb{N}}$. Denote $z_k^{V^\perp} \stackrel{\text{def}}{=} \mathrm{P}_{V^\perp}(z_k)$ the projection of $z_k$ onto $V^\perp$. Moreover, define $\phi_k \stackrel{\text{def}}{=} \frac{1}{2\gamma_k}(\|z_k^{V^\perp} + \gamma_k v^\star\|^2 + \|x_k - x^\star\|^2)$ and the following two auxiliary quantities

$$\xi_k \stackrel{\text{def}}{=} \frac{|\underline{\gamma} - \beta_V|}{2\underline{\gamma}\beta_V}\|z_k - z_{k-1}\|^2 \quad \text{and} \quad \zeta_k \stackrel{\text{def}}{=} \frac{|\gamma_k - \gamma_{k-1}|}{2\underline{\gamma}^2}\|z_k - x^\star\|^2.$$

**Lemma 3.4.** *Considering the non-stationary FDR iteration in (1.6). Suppose that Assumptions (**A.1**)-(**A.5**) hold with $\lambda_k \equiv 1$. Then,*

(i) *We have that $\mathcal{D}_\Phi^{v^\star}(x^\star) = 0$, and*

$$\mathcal{D}_\Phi^{v^\star}(y) \geq 0, \quad \forall y \in \mathbb{R}^n.$$

*Moreover, if $y$ is a solution then $\mathcal{D}_\Phi^{v^\star}(y) = 0$. On the other hand, if $y$ is feasible ($y \in V$) and $\mathcal{D}_\Phi^{v^\star}(y) = 0$, then $y$ is solution.*

(ii) *For the sequence $\{u_k\}_{k\in\mathbb{N}}$, if $v^\star$ is bounded we have*

$$\mathcal{D}_\Phi^{v^\star}(u_{k+1}) + \phi_{k+1} \leq \phi_k + \tfrac{\gamma_{k+1} - \gamma_k}{2}\|v^\star\|^2 + \xi_{k+1} + \zeta_{k+1} < +\infty. \tag{3.5}$$

The proof of the proposition can be found in Section A.

**Remark 3.5.** If we restrict $\gamma_k$ in $]0, \beta_V]$, then the term $\xi_k$ can be discarded in (3.5). If we impose monotonicity assumption on the sequence $\{\gamma_k\}_{k\in\mathbb{N}}$, the term $\zeta_k$ also disappears.

With the above property of $\mathcal{D}_\Phi^{v^\star}(u_k)$, we are able to present the main result on the convergence rate of the Bregman divergence.

**Theorem 3.6.** *Consider the non-stationary FDR iteration (1.6). Suppose that Assumptions (**A.1**)-(**A.5**) hold with $\lambda_k \equiv 1$. If moreover $v^\star$ is bounded, then for any $k \geq 0$,*

$$\inf_{0\leq i\leq k} \mathcal{D}_\Phi^{v^\star}(u_i) = o((k+1)^{-1}) \ \textit{and} \ D(\bar{u}_k) = O(1/(k+1)) \quad \textit{where} \quad \bar{u}_k = \tfrac{1}{k+1}\textstyle\sum_{i=0}^k u_i. \tag{3.6}$$

See Section A for the proof.

**Remark 3.7.**

(i) A typical situation that ensures boundedness of $v^\star$ is when $\partial R(x^\star)$ is bounded. Such requirement can be removed if we choose more carefully the element $v^\star$: for instance we can set $v^\star \stackrel{\text{def}}{=} \frac{x^\star - z^\star}{\gamma}$, which is a sub-gradient of $\nabla F(x^\star) + \partial R(x^\star)$.

(ii) The main difficulty in establishing the convergence rate directly on $\mathcal{D}_\Phi^{v^\star}(u_k)$ (rather that on the best iterate) is that, for $V \subsetneq \mathcal{H}$, we have no theoretical guarantee that $\{\mathcal{D}_\Phi^{v^\star}(u_k)\}_{k\in\mathbb{N}}$ is decreasing, i.e. no descent property on $\mathcal{D}_\Phi^{v^\star}(u_k)$.

## 3.3 Application to FB splitting

Assume now that $V = \mathcal{H}$, in which case problem (1.4) simplifies to

$$\min_{x\in\mathbb{R}^n} \big\{\Phi(x) \stackrel{\text{def}}{=} F(x) + R(x)\big\}.$$

In this case, the FDR iteration (1.6) is nothing but the FB splitting scheme (1.2). The non-relaxed and non-stationary version of it reads as

$$x_{k+1} = \operatorname{prox}_{\gamma_k R}\big(x_k - \gamma_k \nabla F(x_k)\big). \tag{3.7}$$

Specializing to this case the Bregman divergence of (3.4), we get $\mathcal{D}_\Phi^{v^\star}(y) = \Phi(y) - \Phi(x^\star)$, which is simply the objective value error. We have the following result.

**Corollary 3.8.** *Consider the Forward–Backward iteration* (3.7)*. Suppose that conditions* (**A.1**)-(**A.5**) *hold with* $V = \mathcal{H}$ *and* $\lambda_k \equiv 1$*. Then*

$$\Phi(x_k) - \Phi(x^\star) = o(1/(k+1)).$$

See Section A for the proof.

**Remark 3.9.**

(i) The $o(1/(k+1))$ convergence rate for the large choice $\gamma_k \in ]0, 2\beta[$ appears to be new for the FB splitting method.

(ii) For the global convergence of the sequence $\{x_k\}_{k\in\mathbb{N}}$ generated by the non-stationary FB iteration, neither convergence of $\gamma_k$ to $\gamma$ nor summability of $\{|\gamma_k - \gamma|\}_{k\in\mathbb{N}}$ is required. See [14, Theorem 3.4].

# 4 Local linear convergence

From now on, we adopt a finite-dimensional setting where $\mathcal{H} = \mathbb{R}^n$. We investigate the local convergence behaviour of the FDR algorithm. In the sequel, we denote $z^\star \in \operatorname{fix}(\mathscr{F}_\gamma)$ a fixed point of iteration (1.6) and $x^\star = \mathrm{P}_V(z^\star) \in \operatorname{Argmin}(\Phi_V)$ a global minimizer of problem (1.4).

## 4.1 Finite activity identification

We start with the finite activity identification result. Under the condition of Theorem 3.2, we know that $\gamma_k \to \gamma$, $z_k \to z^\star$ and $u_k, x_k \to x^\star$. Moreover, we have the following optimality conditions

$$\frac{x^\star - z^\star}{\gamma} \in \nabla G(x^\star) + \partial R(x^\star) \ \text{ and } \ \frac{z^\star - x^\star}{\gamma} \in V^\perp, x^\star \in V. \tag{4.1}$$

The condition needed for identification result is built upon these monotone inclusions. Since $x_k$ is the projection of $z_k$ onto $V$, we have $x_k \in V$ for all $k \geq 0$. Therefore, we only need to discuss the identification property of $u_k$.

**Theorem 4.1.** *Consider the non-stationary FDR iteration (1.6). Suppose that Assumptions (**A.1**)-(**A.5**) hold, so that $(u_k, x_k, z_k) \to (x^\star, x^\star, z^\star)$ where $z^\star \in \mathrm{fix}(\mathscr{F}_\gamma)$ and $x^\star = \mathrm{P}_V(z^\star) \in \mathrm{Argmin}(\Phi_V)$. Moreover, assume that $R \in \mathrm{PSF}_{x^\star}(\mathcal{M}^R_{x^\star})$ and that the following non-degeneracy condition holds*

$$\frac{x^\star - z^\star}{\gamma} - \nabla G(x^\star) \in \mathrm{ri}\big(\partial R(x^\star)\big). \tag{4.2}$$

*Then,*

(i) *There exists $K \in \mathbb{N}$ such that, for all $k \geq K$, we have $u_k \in \mathcal{M}^R_{x^\star}$.*

(ii) *Moreover,*

(a) *if $\mathcal{M}^R_{x^\star} = x^\star + T^R_{x^\star}$, then $T^R_{u_k} = T^R_{x^\star}$ for every $k \geq K$.*

(b) *If $R$ is locally polyhedral around $x^\star$, then, for every $k \geq K$, $x_k \in \mathcal{M}^R_{x^\star} = x^\star + T^R_{x^\star}$, $T^R_{u_k} = T^R_{x^\star}$, $\nabla_{\mathcal{M}^R_{x^\star}} R(u_k) = \nabla_{\mathcal{M}^R_{x^\star}} R(x^\star)$, and $\nabla^2_{\mathcal{M}^R_{x^\star}} R(u_k) = 0$.*

See Section B.2 for the proof.

**Remark 4.2.** As we mentioned before, for the global convergence of the sequence, approximation errors can be allowed, *i.e.* the proximity operators of $R, J$ and the gradient of $G$ can be computed approximately. However, for the finite activity, we have no identification guarantees for $(u_k, x_k)$ if such an approximation is allowed. For example, if we have $x_k = \mathrm{P}_V(z_k) + \varepsilon_k$ where $\varepsilon_k \in \mathbb{R}^n$ is the error of approximating $\mathrm{P}_V(z_k)$. Then, unless $\varepsilon_k \in V$, we can no longer guarantee that $x_k \in V$.

## 4.2 Locally linearized iteration

With the finite identification result, in the next we show that the globally non-linear fixed-point iteration (3.3) can be locally linearized along the identified manifold $\mathcal{M}^R_{x^\star}$.

Define the following function

$$\overline{R}(u) \stackrel{\text{def}}{=} \gamma R(u) - \langle u,\, x^\star - z^\star - \gamma \nabla G(x^\star)\rangle. \tag{4.3}$$

We have the following key property of $\overline{R}$.

**Lemma 4.3.** *Let $x^\star \in \mathrm{Argmin}(\Phi_V)$, and suppose that $R \in \mathrm{PSF}_{x^\star}(\mathcal{M}^R_{x^\star})$. Then the Riemannian Hessian of $\overline{R}$ at $x^\star$ reads as*

$$H_{\overline{R}} \stackrel{\text{def}}{=} \mathrm{P}_{T^R_{x^\star}} \nabla^2_{\mathcal{M}^R_{x^\star}} \overline{R}(x^\star) \mathrm{P}_{T^R_{x^\star}}, \tag{4.4}$$

*which is symmetric positive semi-definite under either one of the two circumstances:*

(i) *condition (4.2) holds.*

(ii) *$\mathcal{M}^R_{x^\star}$ is an affine subspace.*

*In turn, the matrix $W_{\overline{R}} \stackrel{\text{def}}{=} (\mathrm{Id} + H_{\overline{R}})^{-1}$ is firmly non-expansive.*

**Proof.** Claims (i) and (ii) follow from [30, Lemma 4.3] since $R \in \mathrm{PSF}_{x^\star}(\mathcal{M}^R_{x^\star})$. Consequently, $W_{\overline{R}}$ is symmetric positive definite with eigenvalues in $]0, 1]$. Thus, by virtue of [4, Corollary 4.3(ii)], it is firmly non-expansive. □

From now on, we assume that $F$ (and hence $G$) is locally $C^2$-smooth around $x^\star$. Define $H_G \stackrel{\text{def}}{=} \mathrm{P}_V \nabla^2 F(x^\star) \mathrm{P}_V$, $M_{\overline{R}} \stackrel{\text{def}}{=} \mathrm{P}_{T^R_{x^\star}} W_{\overline{R}} \mathrm{P}_{T^R_{x^\star}}$ and $\mathcal{R}_{M_{\overline{R}}} \stackrel{\text{def}}{=} 2M_{\overline{R}} - \mathrm{Id}$ and the matrices

$$\begin{aligned} \mathscr{M}_\gamma &= \mathrm{Id} + 2M_{\overline{R}}\mathrm{P}_V - M_{\overline{R}} - \mathrm{P}_V - \gamma M_{\overline{R}} H_G = \tfrac{1}{2}\big(\mathcal{R}_{M_{\overline{R}}}\mathcal{R}_V + \mathrm{Id}\big)(\mathrm{Id} - \gamma H_G), \\ \mathscr{M}_{\gamma,\lambda} &= (1-\lambda)\mathrm{Id} + \lambda \mathscr{M}_\gamma. \end{aligned} \tag{4.5}$$

We have the following theorem for the linearized fixed-point formulation of (1.6).

**Theorem 4.4 (Locally linearized iteration).** *Consider the non-stationary FDR iteration (1.6) and suppose that (A.1)-(A.5) hold. If moreover, $\lambda_k \to \lambda \in ]0, \frac{4\beta_V-\gamma}{2\beta_V}[$ and $F$ is locally $C^2$ around $x^\star$, then for all $k$ large enough we have*

$$z_{k+1} - z^\star = \mathcal{M}_{\gamma,\lambda}(z_k - z^\star) + \psi_k + \chi_k, \tag{4.6}$$

*where $\psi_k \stackrel{\text{def}}{=} o(\|z_k - z^\star\|)$ and $\chi_k \stackrel{\text{def}}{=} O(\lambda_k|\gamma_k - \gamma|)$. Both $\psi_k$ and $\chi_k$ vanish when $R$ is locally polyhedral around $x^\star$, $F$ is quadratic and $(\gamma_k, \lambda_k)$ are chosen constants in $]0, 2\beta_V[\times]0, \frac{4\beta_V-\gamma}{2\beta_V}[$.*

See Section B.2 for the proof. Before presenting the local linear convergence result, we need to study the spectral properties of $\mathcal{M}_{\gamma,\lambda}$, which is presented in the lemma below.

**Lemma 4.5 (Convergence properties of $\mathcal{M}_\gamma$).** *Given $\gamma \in ]0, 2\beta_V[$ and $\lambda \in ]0, \frac{4\beta_V-\gamma}{2\beta_V}[$, we have that $\mathcal{M}_\gamma$ is $\frac{2\beta_V}{4\beta_V-\gamma}$-averaged and $\mathcal{M}_{\gamma,\lambda}$ is $\frac{2\beta_V\lambda}{4\beta_V-\gamma}$-averaged. Moreover,*

(i) *$\mathcal{M}_{\gamma,\lambda}$ converges to some matrix $\mathcal{M}_\gamma^\infty$ and, for every $k \in \mathbb{N}$,*

$$\mathcal{M}_{\gamma,\lambda}^k - \mathcal{M}_\gamma^\infty = (\mathcal{M}_{\gamma,\lambda} - \mathcal{M}_\gamma^\infty)^k \quad \textit{and} \quad \rho(\mathcal{M}_{\gamma,\lambda} - \mathcal{M}_\gamma^\infty) < 1.$$

(ii) *Given any $\rho \in ]\rho(\mathcal{M}_{\gamma,\lambda} - \mathcal{M}_\gamma^\infty), 1[$, there is $K$ large enough such that for all $k \geq K$,*

$$\|\mathcal{M}_{\gamma,\lambda}^k - \mathcal{M}_\gamma^\infty\| = O(\rho^k). \tag{4.7}$$

Owing to the convergence property of $\mathcal{M}_{\gamma,\lambda}$, we can further simplify the linearized iteration (4.6).

**Corollary 4.6.** *Consider the non-stationary FDR iteration (1.6) and suppose that it is run under the assumptions of Theorem 4.4. Then the following holds:*

(i) *Iteration (4.6) is equivalent to*

$$(\mathrm{Id} - \mathcal{M}_\gamma^\infty)(z_{k+1} - z^\star) = (\mathcal{M}_{\gamma,\lambda} - \mathcal{M}_\gamma^\infty)(\mathrm{Id} - \mathcal{M}_\gamma^\infty)(z_k - z^\star) + (\mathrm{Id} - \mathcal{M}_\gamma^\infty)\psi_k + \chi_k. \tag{4.8}$$

(ii) *If moreover $R$ is locally polyhedral around $x^\star$, $F$ is quadratic and parameters $(\gamma_k, \lambda_k)$ are chosen constants from $]0, 2\beta_V[\times]0, \frac{4\beta_V-\gamma}{2\beta_V}[$, then*

$$z_{k+1} - z^\star = (\mathcal{M}_{\gamma,\lambda} - \mathcal{M}_\gamma^\infty)(z_k - z^\star). \tag{4.9}$$

See again Section B.2 for the proof.

## 4.3 Local linear convergence

We are now in position to claim local linear convergence of the FDR iterates.

**Theorem 4.7.** *Consider the non-stationary FDR iteration (1.6) and suppose it is run under the conditions of Theorem 4.4. Let be $\rho \in ]\rho(\mathcal{M}_{\gamma,\lambda} - \mathcal{M}_\gamma^\infty), 1[$ and $K \in \mathbb{N}$ such that, for all $k \geq K$, $\|\mathcal{M}_{\gamma,\lambda}^k - \mathcal{M}_\gamma^\infty\| = O(\rho^k)$ (see Lemma 4.5). Then the following holds:*

(i) *If there exists $\eta \in ]0, \rho[$ such that $\lambda_k|\gamma_k - \gamma| = O(\eta^{k-K})$, then*

$$\|(\mathrm{Id} - \mathcal{M}_\gamma^\infty)(z_k - z^\star)\| = O(\rho^{k-K}). \tag{4.10}$$

(ii) *If moreover $R$ is locally polyhedral around $x^\star$, $F$ is quadratic, and*

$$(\gamma_k, \lambda_k) \equiv (\gamma, \lambda) \in ]0, 2\beta_V[\times]0, \frac{4\beta_V - \gamma}{2\beta_V}[,$$

*then we have*

$$\|z_k - z^\star\| \leq \rho^{k-K}\|z_K - z^\star\|. \tag{4.11}$$

**Remark 4.8.**

(i) For the first case of Theorem 4.7, if $\mathscr{M}_\gamma^\infty = 0$ then we obtain the convergence rate directly on $\|z_k - z^\star\|$. Moreover, we can further derive the convergence rate of $\|x_k - x^\star\|$ and $\|u_k - x^\star\|$.

(ii) The condition on $\lambda_k|\gamma_k - \gamma|$ in Theorem 4.7(i) amounts to saying that $\{\gamma_k\}_{k\in\mathbb{N}}$ should converge fast enough to $\gamma$. Otherwise, the local convergence rate would be dominated by that of $\lambda_k|\gamma_k - \gamma|$. Especially, if $\lambda_k|\gamma_k - \gamma|$ converges sub-linearly to 0, then the local convergence rate will eventually become sub-linear. See Figure 2 in the experiments section for a numerical illustration.

(iii) The above result can be easily extended to the case of GFB method, for the sake of simplicity we shall skip the details here. Nevertheless, numerical illustrations will be provided in Section 5.

### 4.4 Extension to the three-operator splitting

So far, we have presented the global and local convergence analysis of the FDR method. As we recalled in the introduction, FDR method is closely related to the three-operator splitting method (TOS) [17]. Therefore, it would be interesting to extend the obtained result to TOS method. However, extending the global convergence result to TOS is far from straightforward. Hence, in the following, we mainly focus on the local linear convergence results.

For the sake of notational simplicity, we rewrite problem (1.9) as

$$\min_{x\in\mathbb{R}^n} \big\{\Psi(x) = F(x) + R(x) + J(x)\big\}, \tag{4.12}$$

where we suppose the following assumptions:

**(B.1)** $J, R \in \Gamma_0(\mathbb{R}^n)$.

**(B.2)** $F : \mathcal{H} \to \mathbb{R}$ is convex continuously differentiable with $\nabla F$ being $(1/\beta)$-Lipschitz continuous.

**(B.3)** $\mathrm{Argmin}(\Psi) \neq \emptyset$, *i.e.* the set of minimizers is not empty.

Correspondingly, the TOS iteration (1.10) becomes

$$\begin{cases} u_{k+1} = \mathrm{prox}_{\gamma R}\big(2x_k - z_k - \gamma\nabla F(x_k)\big) \\ z_{k+1} = z_k + \lambda_k(u_{k+1} - x_k), \\ x_{k+1} = \mathrm{prox}_{\gamma J}(z_{k+1}). \end{cases} \tag{4.13}$$

We suppose the following assumption on the algorithm parameters:

**(B.4)** The (constant) step-size verifies $\gamma \in ]0, 2\beta[$ and the sequence of relaxation parameters $\{\lambda_k\}_{k\in\mathbb{N}}$ is such that $\sum_{k\in\mathbb{N}} \lambda_k(\frac{4\beta-\gamma}{2\beta} - \lambda_k) = +\infty$.

The fixed-point operator of TOS reads as

$$\mathscr{T}_\gamma = \mathrm{Id} - \mathrm{prox}_{\gamma R} + \mathrm{prox}_{\gamma J}(2\mathrm{prox}_{\gamma R} - \mathrm{Id} - \gamma\nabla F \circ \mathrm{prox}_{\gamma J}) \quad \text{and} \quad \mathscr{T}_{\gamma,\lambda_k} = (1-\lambda_k)\mathrm{Id} + \lambda_k\mathscr{T}_\gamma. \tag{4.14}$$

Differently from $\mathscr{F}_\gamma$ (see (3.1)), $\mathscr{T}_\gamma$ cannot be simplified into a compact form. We have the following lemma for the convergence of the TOS method.

**Lemma 4.9 (Global convergence of TOS).** *Consider the TOS iteration (4.13) and the fixed-point operator (4.14). Suppose that Assumptions (**B.1**)-(**B.4**) hold. Then*

(i) *The operator $\mathscr{T}_\gamma$ is $\frac{2\beta}{4\beta-\gamma}$-averaged non-expansive.*

(ii) *The sequence $z_k$ converges to some $z^\star$ in $\mathrm{fix}(\mathscr{T}_\gamma)$; moreover, both $u_k$ and $x_k$ converge to $x^\star \overset{\text{def}}{=} \mathrm{prox}_{\gamma J}(z^\star)$, that is a global minimizer of $\mathrm{Argmin}(\Psi)$.*

**Proof.** See Proposition 2.1 and Theorem 2.1 of [17]. □

Similar to (4.1), under Lemma 4.9, we have the following optimality condition at convergence:

$$\frac{x^\star - z^\star}{\gamma} \in \nabla F(x^\star) + \partial R(x^\star) \ \text{ and } \ \frac{z^\star - x^\star}{\gamma} \in \partial J(x^\star). \tag{4.15}$$

Following the footprint of Section 4.1-4.3, in the next we present the local linear convergence of TOS method.

**Finite activity identification** We start with the finite identification result, since $J$ is no longer the indicator function of a subspace, we have the identification result for both $u_k$ and $x_k$.

**Corollary 4.10.** *Consider the TOS iteration (4.13). Suppose it is run under the Assumptions (**B.1**)-(**B.4**), so that $(u_k, x_k, z_k) \to (x^\star, x^\star, z^\star)$ where $z^\star \in \mathrm{fix}(\mathscr{T}_\gamma)$ and $x^\star \stackrel{\text{def}}{=} \mathrm{prox}_{\gamma J}(z^\star) \in \mathrm{Argmin}(\Psi)$. Moreover, let $R \in \mathrm{PSF}_{x^\star}(\mathcal{M}^R_{x^\star})$, $J \in \mathrm{PSF}_{x^\star}(\mathcal{M}^J_{x^\star})$, and the following non-degeneracy condition holds*

$$\frac{x^\star - z^\star}{\gamma} - \nabla F(x^\star) \in \mathrm{ri}\big(\partial R(x^\star)\big) \ \text{ and } \ \frac{z^\star - x^\star}{\gamma} \in \mathrm{ri}\big(\partial J(x^\star)\big). \tag{4.16}$$

*Then, there exists $K \in \mathbb{N}$ such that $(u_k, x_k) \in \mathcal{M}^R_{x^\star} \times \mathcal{M}^J_{x^\star}$ for every $k \geq K$.*

**Local linearized iteration** The next step is to linearize the TOS iteration. Define functions

$$\widetilde{J}(x) \stackrel{\text{def}}{=} \gamma J(x) - \langle x,\, z^\star - x^\star \rangle \ \text{ and } \ \widetilde{R}(u) \stackrel{\text{def}}{=} \gamma R(u) - \langle u,\, x^\star - z^\star - \gamma \nabla F(x^\star) \rangle. \tag{4.17}$$

We start with the following result, corollary from Lemma 4.3.

**Corollary 4.11.** *Suppose that $J \in \mathrm{PSF}_{x^\star}(\mathcal{M}^J_{x^\star})$ and $R \in \mathrm{PSF}_{x^\star}(\mathcal{M}^R_{x^\star})$. Then their Riemannian Hessians at $x^\star$ read*

$$H_{\widetilde{J}} \stackrel{\text{def}}{=} \mathrm{P}_{T^J_{x^\star}} \nabla^2_{\mathcal{M}^J_{x^\star}} \widetilde{J}(x^\star) \mathrm{P}_{T^J_{x^\star}} \ \text{ and } \ H_{\widetilde{R}} \stackrel{\text{def}}{=} \mathrm{P}_{T^R_{x^\star}} \nabla^2_{\mathcal{M}^R_{x^\star}} \widetilde{R}(x^\star) \mathrm{P}_{T^R_{x^\star}}, \tag{4.18}$$

*which are symmetric positive semi-definite under either of the following circumstances:*

(i) *condition (4.16) holds.*

(ii) *$\mathcal{M}^J_{x^\star}$ and $\mathcal{M}^R_{x^\star}$ are affine subspaces.*

*In turn, the matrices*

$$W_{\widetilde{J}} \stackrel{\text{def}}{=} (\mathrm{Id} + H_{\widetilde{J}})^{-1} \ \text{ and } \ W_{\widetilde{R}} \stackrel{\text{def}}{=} (\mathrm{Id} + H_{\widetilde{R}})^{-1} \tag{4.19}$$

*are both firmly non-expansive.*

Now assume $F$ is locally $C^2$-smooth around $x^\star$, and define $H_F \stackrel{\text{def}}{=} \nabla^2 F(x^\star)$. Define also $M_{\widetilde{J}} \stackrel{\text{def}}{=} \mathrm{P}_{T^J_{x^\star}} W_{\widetilde{J}} \mathrm{P}_{T^J_{x^\star}}$ and $M_{\widetilde{R}} \stackrel{\text{def}}{=} \mathrm{P}_{T^R_{x^\star}} W_{\widetilde{R}} \mathrm{P}_{T^R_{x^\star}}$, and the matrices

$$\mathscr{L}_\gamma = \mathrm{Id} + 2M_{\widetilde{R}} M_{\widetilde{J}} - M_{\widetilde{R}} - M_{\widetilde{J}} - \gamma M_{\widetilde{R}} H_F M_{\widetilde{J}} \ \text{ and } \ \mathscr{L}_{\gamma,\lambda} = (1-\lambda)\mathrm{Id} + \lambda \mathscr{L}_\gamma.$$

**Lemma 4.12.** *The matrix $\mathscr{L}_\gamma$ is $\frac{2\beta}{4\beta-\gamma}$-averaged non-expansive.*

**Proof.** See [17, Proposition 2.1]. □

The above lemma entails that $\mathscr{L}_\gamma, \mathscr{L}_{\gamma,\lambda}$ are convergent, hence the spectral properties result in Lemma 4.5 applies to them. Denote $\mathscr{L}^\infty_\gamma \stackrel{\text{def}}{=} \lim_{k\to+\infty} \mathscr{L}^k_{\gamma,\lambda}$. Now we are able to present the result on the local linearization of iteration (4.13) along the identified manifolds.

**Corollary 4.13 (Local linearized iteration).** *Consider the TOS iteration (4.13). Suppose it is run under Assumptions (**B.1**)-(**B.4**), that* $\lambda_k \to \lambda \in ]0, \frac{4\beta-\gamma}{2\beta}[$*, and that* $F$ *is locally* $C^2$ *around* $x^\star$*. Then we have*

$$z_{k+1} - z^\star = \mathscr{L}_{\gamma,\lambda}(z_k - z^\star) + o(\|z_k - z^\star\|). \tag{4.20}$$

*If moreover* $J, R$ *are locally polyhedral around* $x^\star$*,* $F$ *is quadratic and* $\lambda_k \equiv \lambda \in ]0, \frac{4\beta-\gamma}{2\beta}[$ *is chosen constant, then the term* $o(\|z_k - z^\star\|)$ *vanishes.*

We can also specialize Corollary 4.6 to this context, however we choose to skip it owing to its obviousness.

**Local linear convergence** Finally, we are able to to present the local linear convergence for (4.13).

**Corollary 4.14.** *Consider the TOS iteration (4.13). Suppose it is run under Assumptions (**B.1**)-(**B.4**), that* $\lambda_k \to \lambda \in ]0, \frac{4\beta-\gamma}{2\beta}[$*, and that* $F$ *is locally* $C^2$ *around* $x^\star$*. Then the following holds:*

(i) *Given any* $\rho \in ]\rho(\mathscr{L}_{\gamma,\lambda} - \mathscr{L}_\gamma^\infty), 1[$*, there exists* $K \in \mathbb{N}$ *large enough such that, for all* $k \geq K$*,*

$$\|(\mathrm{Id} - \mathscr{L}_\gamma^\infty)(z_k - z^\star)\| = O(\rho^{k-K}).$$

(ii) *If moreover* $J, R$ *are locally polyhedral around* $x^\star$*,* $F$ *is quadratic and* $\lambda_k \equiv \lambda \in ]0, \frac{4\beta-\gamma}{2\beta}[$ *is chosen constant, then there exists* $K \in \mathbb{N}$ *such that, for all* $k \geq K$*,*

$$\|z_k - z^\star\| \leq \rho^{k-K}\|z_K - z^\star\|.$$

# 5 Numerical experiments

In this section, we illustrate our theoretical results on concrete examples arising from statistics, and signal/image processing applications.

## 5.1 Examples of partly smooth functions

Table 1 provides some examples of partly smooth functions that we will use throughout this section. More details about them can be found in [30, Section 5] and references therein.

Table 1: Examples of partly smooth functions. For $x \in \mathbb{R}^n$ and some subset of indices $\mathfrak{b} \subset \{1, \ldots, n\}$, $x_\mathfrak{b}$ is the restriction of $x$ to the entries indexed in $\mathfrak{b}$. $D_{\mathrm{DIF}}$ stands for the finite differences operator.

| Function | Expression | Partial smooth manifold |
|---|---|---|
| $\ell_1$-norm | $\|x\|_1 = \sum_{i=1}^n |x_i|$ | $\mathcal{M} = T_x = \{z \in \mathbb{R}^n : I_z \subseteq I_x\}$, $I_x = \{i : x_i \neq 0\}$ |
| $\ell_{1,2}$-norm | $\sum_{i=1}^m \|x_{\mathfrak{b}_i}\|$ | $\mathcal{M} = T_x = \{z \in \mathbb{R}^n : I_z \subseteq I_x\}$, $I_x = \{i : x_{\mathfrak{b}_i} \neq 0\}$ |
| $\ell_\infty$-norm | $\max_{i=\{1,\ldots,n\}} |x_i|$ | $\mathcal{M} = T_x = \{z \in \mathbb{R}^n : z_{I_x} \in \mathbb{R}\mathrm{sign}(x_{I_x})\}$, $I_x = \{i : |x_i| = \|x\|_\infty\}$ |
| TV semi-norm | $\|x\|_{\mathrm{TV}} = \|D_{\mathrm{DIF}}x\|_1$ | $\mathcal{M} = T_x = \{z \in \mathbb{R}^n : I_{D_{\mathrm{DIF}}z} \subseteq I_{D_{\mathrm{DIF}}x}\}$, $I_{D_{\mathrm{DIF}}x} = \{i : (D_{\mathrm{DIF}}x)_i \neq 0\}$ |
| Nuclear norm | $\|x\|_* = \sum_{i=1}^r \sigma(x)$ | $\mathcal{M} = \{z \in \mathbb{R}^{n_1 \times n_2} : \mathrm{rank}(z) = \mathrm{rank}(x) = r\}$, $\sigma(x)$ singular values of $x$ |

The $\ell_1$, $\ell_\infty$-norms and the anisotropic TV semi-norm are all polyhedral functions, hence the corresponding Riemannian Hessians are simply 0. The $\ell_{1,2}$-norm is not polyhedral yet partly smooth relative to a subspace; the nuclear norm is partly smooth relative to the manifold of fixed-rank matrices, which is not anymore flat. The Riemannian Hessian of these two functions are non-trivial and can be computed following [43].

### 5.2 Global convergence rate of the Bregman distance

We first demonstrate, numerically, the global $o(1/k)$ convergence rate of the Bregman divergence of Section 3. Towards this goal, we consider the linearly constrained LASSO problem [22, 20, 23], namely

$$\min_{x\in\mathbb{R}^n} \frac{1}{2}\|\mathcal{K}x - f\|^2 + \mu\|x\|_1 + \iota_V(x). \tag{5.1}$$

In the latter, $\mu > 0$ is a weight parameter, $\mathcal{K} : \mathbb{R}^n \to \mathbb{R}^m$ is a linear operator, and $V \subset \mathbb{R}^n$ is a subspace.

It is immediate to see that problem (5.1) falls within the scope of the FDR algorithm. We set $\lambda_k \equiv 1$, and the step-size $\gamma_k \equiv \frac{1}{7\|\mathcal{K}\|^2}$.

The convergence profile of $\min_{0\le i\le k} \mathcal{D}_\Phi^{v^\star}(u_i)$ is shown in Figure 1(a). The plot is in log-log scale, where the red line corresponds to the sub-linear $O(1/k)$ rate and the black line is $\min_{0\le i\le k} \mathcal{D}_\Phi^{v^\star}(u_i)$. One can then confirm numerically the prediction of Theorem 3.6.

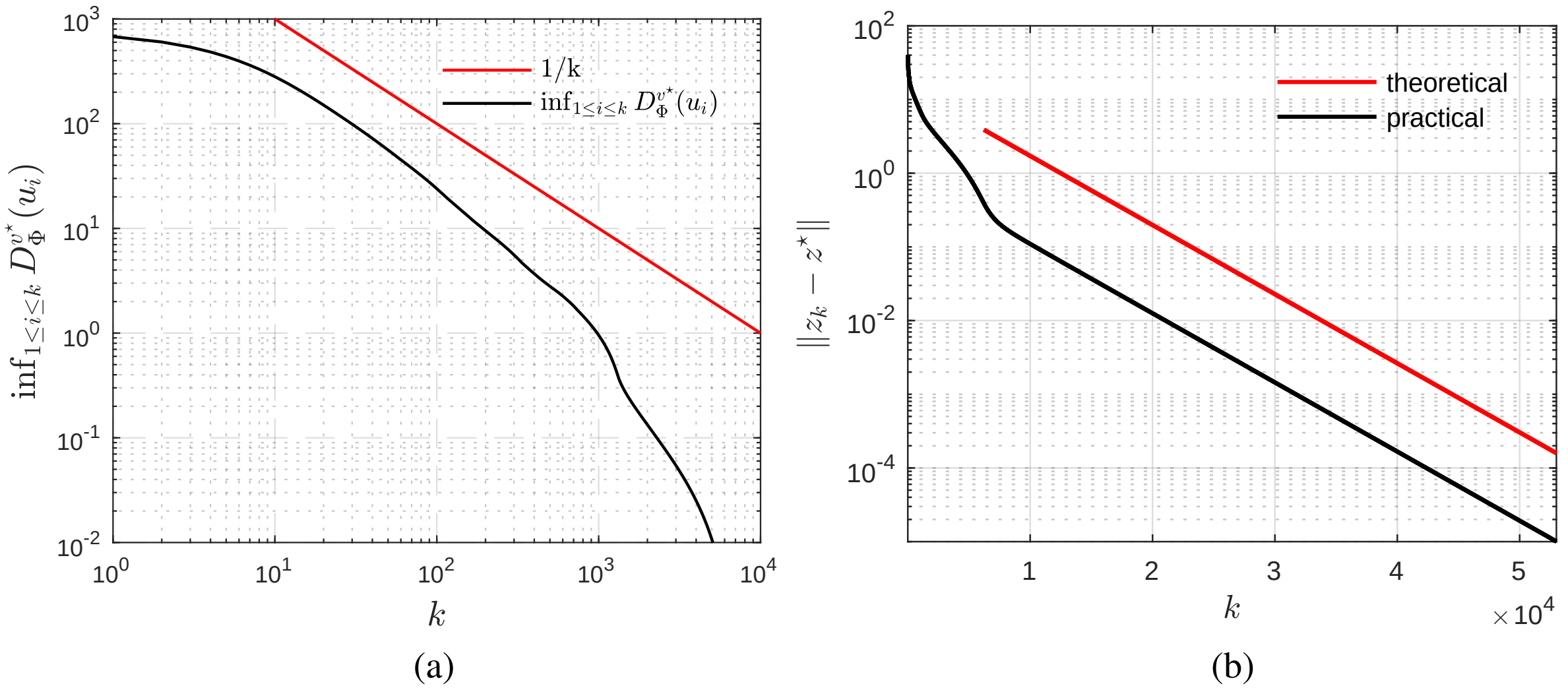

Figure 1: (a): convergence rate of the Bregman distance $\inf_{0\le i\le k} \mathcal{D}_\Phi^{v^\star}(u_i)$ of FDR algorithm for solving the subspace constrained LASSO problem (5.1). (b): convergence rate of the $\|z_k - z^\star\|$ of FDR algorithm for solving the subspace constrained LASSO problem (5.1).

However, it can be observed that beyond some iteration, typically $k > 10^3$ the convergence rate regime changes to become linear. We argue in the following section that this is likely to be due to finite activity identification since the $\ell_1$-norm is partly smooth (in fact even polyhedral) and that, for all $k$ large enough, FDR enters into a local linear convergence regime.

### 5.3 Local linear convergence of FDR

Following the above discussion, in Figure 1(b) we present the local linear convergence of FDR algorithm in terms of $\|z_k - z^\star\|$ as we are in the scope of Theorem 4.7(ii). We use the same parameters setting as in Figure 1(a). The red line stands for the estimated rate (see Theorem 4.7), while the black line is the numerical observation. The starting point of the red line is the number of iteration where $u_k$ identifies the manifold $\mathcal{M}_{x^\star}^R$. As shown in the figure, we indeed have local linear convergence behaviour of $\|z_k - z^\star\|$. Moreover,

since $F = \frac{1}{2}\|\mathcal{K}x - f\|^2$ is quadratic and $R = \|x\|_1$ is polyhedral, our theoretical rate estimation is tight, *i.e.* the red line has the same slope as the black line.

We now investigate numerically the convergence behaviour of the non-stationary version of FDR and compare it to the stationary one. We fix $\lambda_k \equiv 1$, *i.e.* the iteration is unrelaxed. The stationary FDR algorithm is run with $\gamma = \beta$. For the non-stationary ones, four choices of $\gamma_k$ are considered:

$$\text{Case 1: } \gamma_k = (1 + \tfrac{1}{k^{1.1}})\beta, \qquad \text{Case 2: } \gamma_k = (1 + \tfrac{1}{k^2})\beta,$$
$$\text{Case 3: } \gamma_k = (1 + 0.999^k)\beta, \qquad \text{Case 4: } \gamma_k = (1 + 0.5^k)\beta.$$

Obviously, we have $\gamma_k \to \gamma = \beta$ and $\sum_{k\in\mathbb{N}} |\gamma_k - \gamma| < +\infty$ for all the four cases. Problem (5.1) is considered again. The comparison results are displayed in Figure 2(a).

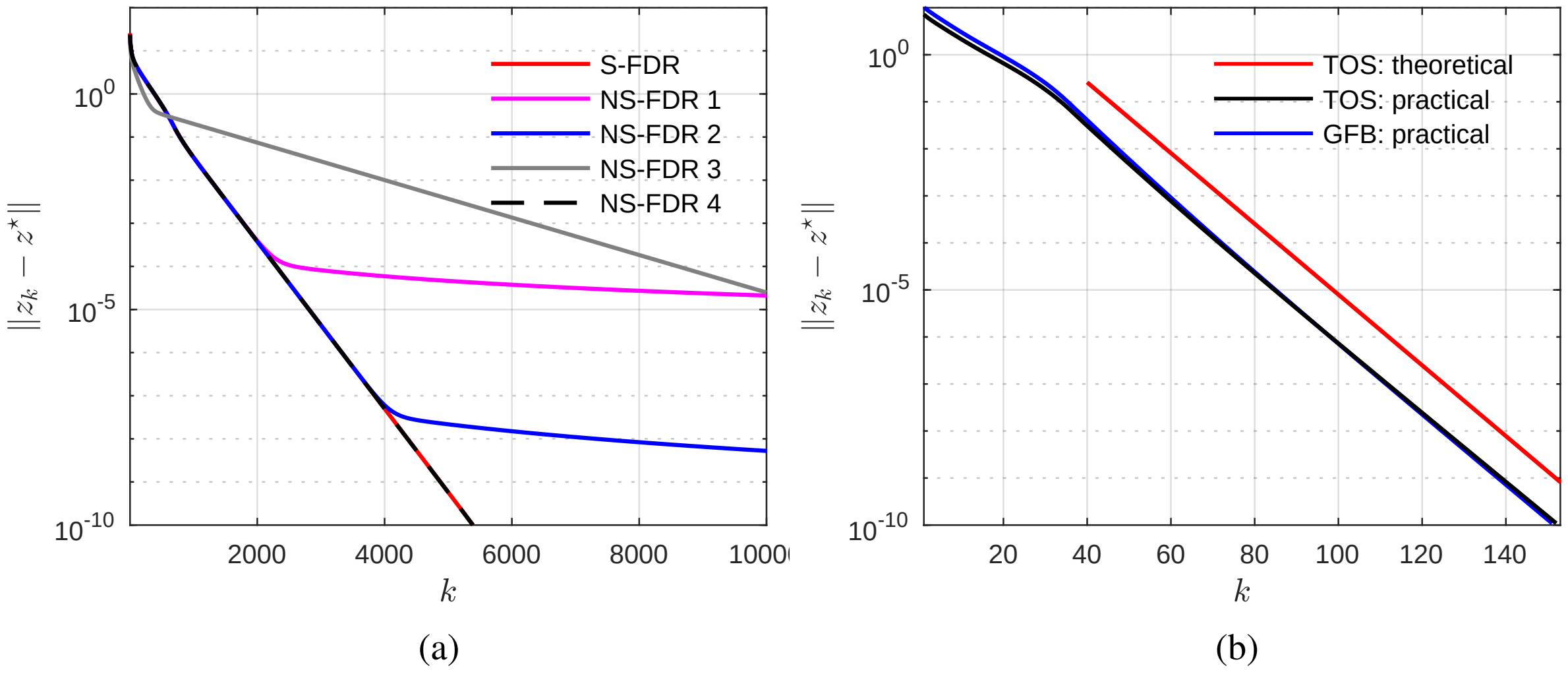

Figure 2: Comparison in termes of the decay of $\|z_k - z^\star\|$ between stationary ("S-FDR") and non-stationary FDR ("NS-FDR X", X stands for Case X).

For the stationary iteration, the local convergence rate of FDR is $\rho = 0.9955$. We can make the following observations from the comparison:

- In agreement with our analysis, the local convergence behaviour of the non-stationary iteration is no better than the stationary one. This contrasts with the global behaviour where non-stationarity could be beneficial (see last comment hereafter);
- As argued in Remark 4.8(ii), the convergence rate is eventually controlled by the error $|\gamma_k - \gamma|$, except for "Case 4", Indeed, 0.5 is strictly smaller than the local linear rate of the stationary version (*i.e.* $|\gamma_k - \gamma| = o(\|z_k - z^\star\|)$);
- The non-stationary FDR seems to lead to faster identification, typically for "Case 3". This is the effect of bigger step-size at the early stage of the algorithm.

## 5.4 Local linear convergence of GFB and TOS

To conclude the numerical experiments, we demonstrate the local convergence behaviour of TOS method. Consider the following problem, recovering a vector which is group sparse and piece-wise constant inside

each non-zero group:

$$\min_{x\in\mathbb{R}^n} \mu_1\|x\|_{1,2} + \frac{1}{2}\|\mathcal{K}x - f\|^2 + \mu_2\|D_{\mathrm{DIF}}x\|_1.$$

Parameters $\mu_1, \mu_2 > 0$ are tradeoff weights. This problem is a special case of (4.12) if we let $F = \frac{1}{2}\|\mathcal{K}\cdot - b\|^2$, $R = \mu_1\|\cdot\|_{1,2}$ and $J = \mu_2\|D_{\mathrm{DIF}}\cdot\|_1$. Hence it can be solved by the TOS scheme (4.13) and also by the GFB algorithm (1.8).

Figure 2(b) shows the local convergence behaviour of the TOS method solving problem (5.4). The red line stands for the theoretical rate (estimated by Corollary 4.14), while the black line is the numerical observation. Similar to the observation of FDR method in Figure 1(b), local linear convergence of TOS occurs and our theoretical rate estimation is rather tight. We also plot the convergence profile of the GFB method (blue line). It can be observed that the performances of GFB and TOS are almost the same.

# 6 Conclusion

In this paper, we studied global and local convergence properties of the Forward–Douglas–Rachford method. Globally, we established an $o(1/k)$ convergence rate of the best iterate and $O(1/k)$ ergodic rate in terms of a Bregman divergence criterion designed for the method. We also specialized the result to the case of Forward–Backward splitting method, for which we have shown that the objective function of the method converges at an $o(1/k)$ rate. Then, locally, we proved the linear convergence of the sequence when the involved functions are moreover partly smooth. In particular, we demonstrated that the method identifies the active manifolds in finite time and that then it converges locally linearly at a rate that we characterized precisely. We also extended the local linear convergence result to the case of Three-Operator-Splitting method. Out numerical experiments supported the theoretical findings.

# Acknowledgements

CM was supported by CONICYT scholarship CONICYT-PCHA/Doctorado Nacional/2016. JL was supported by the European Research Council (ERC project SIGMA-Vision), and Leverhulme Trust project "Breaking the non-convexity barrier", the EPSRC grant "EP/M00483X/1", EPSRC centre "EP/N014588/1", the Cantab Capital Institute for the Mathematics of Information, and the Global Alliance project "Statistical and Mathematical Theory of Imaging". JF was partly supported by Institut Universitaire de France.

# A Proofs of Section 3

**Proof of Theorem 3.2.** To prove the convergence of the non-stationary iteration, we only need to check the conditions listed in [29, Theorem 4].

Owing to Lemma 3.1, given $\gamma_k \in [0, 2\beta_V]$, we have that $\mathscr{F}_{\gamma_k}$ is $\alpha_k$-averaged with $\alpha_k = \frac{2\beta_V}{4\beta_V - \gamma_k}$. This means that there exists a non-expansive operator $\mathscr{R}_{\gamma_k}$ such that $\mathscr{F}_{\gamma_k} = \alpha_k\mathscr{R}_{\gamma_k} + (1-\alpha_k)\mathrm{Id}$. Similarly, for $\gamma \in [0, 2\beta_V]$, we have that $\mathscr{F}_\gamma$ is $\alpha$-averaged with $\alpha = \frac{2\beta_V}{4\beta_V - \gamma}$ and so that there exists a non-expansive operator $\mathscr{R}_\gamma$ such that $\mathscr{F}_\gamma = \alpha\mathscr{R}_\gamma + (1-\alpha)\mathrm{Id}$. Provided $z_k$, define the error term $e_k = (\mathscr{F}_{\gamma_k} - \mathscr{F}_\gamma)(z_k)$. Then iteration (3.3) can be written as

$$\begin{aligned} z_{k+1} = (1-\lambda)z_k + \lambda_k\mathscr{F}_{\gamma_k}(z_k) &= (1-\lambda)z_k + \lambda_k\big(\mathscr{F}_\gamma(z_k) + (\mathscr{F}_{\gamma_k}(z_k) - \mathscr{F}_\gamma(z_k))\big) \\ &= (1-\lambda)z_k + \lambda_k\big(\mathscr{F}_\gamma(z_k) + e_k\big). \end{aligned} \tag{A.1}$$

From the assumptions (**A.1**)-(**A.5**), we can derive the following results:

(i) Recall that $\operatorname{Argmin}(\Phi_V)$ is nothing but $\mathrm{P}_V(\mathrm{fix}(\mathscr{F}_\gamma))$. Then, as (A.4) assumes the set of minimizers of (1.4) is nonempty, so it is the set $\mathrm{fix}(\mathscr{F}_\gamma)$.

(ii) Owing to Lemma 3.1, we have that $\mathscr{F}_{\gamma_k,\lambda_k}$ is $(\alpha_k\lambda_k)$-averaged non-expansive.

(iii) Owing to the averageness of $\mathscr{F}_\gamma$ and $\mathscr{F}_{\gamma_k}$, we have

$$\mathscr{R}_\gamma = \mathrm{Id} + \frac{1}{\alpha}(\mathscr{F}_\gamma - \mathrm{Id}) \ \text{ and } \ \mathscr{R}_{\gamma_k} = \mathrm{Id} + \frac{1}{\alpha_k}(\mathscr{F}_{\gamma_k} - \mathrm{Id}).$$

Let $\rho > 0$ be a positive number. Then, $\forall z \in \mathcal{H}$ such that $\|z\| \le \rho$, we have

$$\begin{aligned}
\|\mathscr{R}_{\gamma_k}(z) - \mathscr{R}_\gamma(z)\| &= \|\frac{1}{\alpha_k}(\mathscr{F}_{\gamma_k} - \mathrm{Id})(z) - \frac{1}{\alpha}(\mathscr{F}_\gamma - \mathrm{Id})(z)\| \\
&= \|\frac{1}{\alpha_k}(\mathscr{F}_{\gamma_k} - \mathrm{Id})(z) - \frac{1}{\alpha_k}(\mathscr{F}_\gamma - \mathrm{Id})(z) + \frac{1}{\alpha_k}(\mathscr{F}_\gamma - \mathrm{Id})(z) - \frac{1}{\alpha}(\mathscr{F}_\gamma - \mathrm{Id})(z)\| \\
&\le |\frac{1}{\alpha_k} - \frac{1}{\alpha}|\|(\mathscr{F}_\gamma - \mathrm{Id})(z)\| + \frac{1}{\alpha_k}\|(\mathscr{F}_{\gamma_k} - \mathrm{Id})(z) - (\mathscr{F}_\gamma - \mathrm{Id})(z)\| \\
&\le \frac{|\gamma_k - \gamma|}{2\beta_V}\big(2\rho + \|\mathscr{F}_\gamma(0)\|\big) + \frac{1}{\alpha_k}\|(\mathscr{F}_{\gamma_k} - \mathscr{F}_\gamma)(z)\|.
\end{aligned} \tag{A.2}$$

Given $\gamma \in ]0, 2\beta_V[$, define the following two operators

$$\mathscr{F}_{1,\gamma} = \frac{1}{2}(\mathrm{Id} + \mathcal{R}_{\gamma R} \circ \mathcal{R}_V) \ \text{ and } \ \mathscr{F}_{2,\gamma} = \mathrm{Id} - \gamma\nabla G.$$

Then $\mathscr{F}_{1,\gamma}$ is firmly expansive and $\mathscr{F}_{2,\gamma}$ is $\frac{\gamma}{2\beta_V}$-averaged (Lemma 2.7). Now we have

$$\begin{aligned}
\|\mathscr{F}_{\gamma_k}(z) - \mathscr{F}_\gamma(z)\| &= \|\mathscr{F}_{1,\gamma_k}\mathscr{F}_{2,\gamma_k}(z) - \mathscr{F}_{1,\gamma_k}\mathscr{F}_{2,\gamma}(z) + \mathscr{F}_{1,\gamma_k}\mathscr{F}_{2,\gamma}(z) - \mathscr{F}_{1,\gamma}\mathscr{F}_{2,\gamma}(z)\| \\
&\le \|\mathscr{F}_{1,\gamma_k}\mathscr{F}_{2,\gamma_k}(z) - \mathscr{F}_{1,\gamma_k}\mathscr{F}_{2,\gamma}(z)\| + \|\mathscr{F}_{1,\gamma_k}\mathscr{F}_{2,\gamma}(z) - \mathscr{F}_{1,\gamma}\mathscr{F}_{2,\gamma}(z)\| \\
&\le \underbrace{\|\mathscr{F}_{2,\gamma_k}(z) - \mathscr{F}_{2,\gamma}(z)\|}_{\text{Term 1}} + \underbrace{\|\mathscr{F}_{1,\gamma_k}\mathscr{F}_{2,\gamma}(z) - \mathscr{F}_{1,\gamma}\mathscr{F}_{2,\gamma}(z)\|}_{\text{Term 2}}.
\end{aligned} \tag{A.3}$$

For the first term of (A.3),

$$\begin{aligned}
\|\mathscr{F}_{2,\gamma_k}(z) - \mathscr{F}_{2,\gamma}(z)\| &= |\gamma_k - \gamma|\|\nabla G(z)\| \\
\text{(Triangle inequality and } \nabla G \text{ is } \beta_V^{-1}\text{-Lipschitz)} &\le |\gamma_k - \gamma|(\beta_V^{-1}\rho + \|\nabla G(0)\|),
\end{aligned} \tag{A.4}$$

where $\nabla G(0)$ is obviously bounded. Now consider the second term of (A.3). Denoting $z^V = \mathrm{P}_V(z)$ and $z^{V^\perp} = z - z^V$, it can be derived that

$$v = \mathscr{F}_{1,\gamma}\mathscr{F}_{2,\gamma}(z) \iff v = z^{V^\perp} + \mathrm{prox}_{\gamma R}(z^V - z^{V^\perp} - \gamma\nabla G(z^V)).$$

Denote $y = z^V - z^{V^\perp} - \gamma\nabla G(z^V)$. Then we have

$$\mathscr{F}_{1,\gamma_k}\mathscr{F}_{2,\gamma}(z) - \mathscr{F}_{1,\gamma}\mathscr{F}_{2,\gamma}(z) = \mathrm{prox}_{\gamma_k R}(y) - \mathrm{prox}_{\gamma R}(y).$$

Denote $w_k = \mathrm{prox}_{\gamma_k R}(y)$ and $w = \mathrm{prox}_{\gamma R}(y)$. Using the classical resolvent equation [9] and non-expansiveness of the proximity operator, we have

$$\|w_k - w\| = \|\mathrm{prox}_{\gamma_k R}(\gamma_k/\gamma y + (1 - \gamma_k/\gamma)w) - \mathrm{prox}_{\gamma_k R}(y)\| \le \|(1 - \gamma_k/\gamma)(y - w)\|.$$

This, together with firm non-expansiveness of $\mathrm{Id} - \mathrm{prox}_{\gamma R}$ yields

$$\|w - w_k\| \le \frac{|\gamma_k - \gamma|}{\underline{\gamma}}\|y - w\| = \frac{|\gamma_k - \gamma|}{\underline{\gamma}}\|(\mathrm{Id} - \mathrm{prox}_{\gamma R})y\| \le \frac{|\gamma_k - \gamma|}{\underline{\gamma}}(\|y\| + \|\mathrm{prox}_{\gamma R}(0)\|). \tag{A.5}$$

Using the triangle inequality, the Pythagorean theorem and non-expansiveness of $\beta_V\nabla G$, we obtain

$$\begin{aligned}
\|y\| \le \|z^V - z^{V^\perp}\| + \gamma\|\nabla G(z^V)\| &\le \rho + \gamma\|\nabla G(z^V) - \nabla G(0)\| + \gamma\|\nabla G(0)\| \\
&\le \rho + \gamma\beta_V^{-1}\|z\| + \gamma\|B(0)\| \le \rho + \bar{\gamma}\beta_V^{-1}\rho + \bar{\gamma}\|B(0)\|.
\end{aligned} \tag{A.6}$$

Define $\Delta_{k,\rho} \stackrel{\text{def}}{=} \sup_{\|z\|\leq\rho} \|\mathscr{R}_{\gamma_k}(z) - \mathscr{R}_\gamma(z)\|$. Then, putting together (A.2), (A.4), (A.5) and (A.6), we get that $\forall \rho \in [0, +\infty[$

$$\textstyle\sum_{k\in\mathbb{N}}\lambda_k\alpha_k\Delta_{k,\rho} = \sum_{k\in\mathbb{N}}\lambda_k\alpha_k \sup\limits_{\|z\|\leq\rho} \|\mathscr{R}_{\gamma_k}(z) - \mathscr{R}_\gamma(z)\| \leq C \sum_{k\in\mathbb{N}} \lambda_k|\gamma_k - \gamma| < +\infty,$$

where $C = \frac{2\rho+\|\mathscr{F}_\gamma(0)\|}{4\beta_V - \bar{\gamma}} + \frac{\rho}{\beta_V}(1 + \frac{\beta_V}{\underline{\gamma}} + \frac{\bar{\gamma}}{\underline{\gamma}}) + (1 + \frac{\bar{\gamma}}{\underline{\gamma}})\|\nabla G(0)\| + \frac{1}{\underline{\gamma}}\|\mathrm{prox}_{\gamma R}(0)\| < +\infty$.
As we verified that all the conditions listed in [29, Theorem 4] are met for the non-stationary FDR iteration, weak convergence of the sequence $\{z_k\}_{k\in\mathbb{N}}$ then follows. In turn, since $\mathrm{P}_V$ is linear, weak convergence of $\{x_k\}_{k\in\mathbb{N}}$ is obtained.

For the sequence $\{u_k\}_{k\in\mathbb{N}}$, observe from the second equation in (1.6) that $u_{k+1} = (z_{k+1} - z_k)\lambda_k + x_k$, and thus $\|u_{k+1} - x_k\| \leq (\|z_{k+1} - z_k\|)/\lambda_k$. It follows from $\|z_{k+1} - z_k\| \to 0$ and the condition $\inf_{k\in\mathbb{N}} \lambda_k > 0$ that $u_{k+1} - x_k$ converges strongly to 0. We thus obtain weak convergence of $u_k$ to the same weak cluster point as $x_k$.

If $\mathcal{H}$ is finite-dimensional, we use the optimality conditions (4.1) and argue as for (A.5) to obtain

$$\begin{aligned} \|u_{k+1} - x^\star\| &\leq \frac{|\gamma_k - \gamma|}{\underline{\gamma}}(\|2(x_k - x^\star) + (z_k - z^\star) + \gamma(\nabla G(x_k) - \nabla G(x^\star))\| + \|\mathrm{prox}_{\gamma R}(0)\|) \\ &\leq \frac{|\gamma_k - \gamma|}{\underline{\gamma}}\big((2 + \bar{\gamma}\beta_V)\|x_k - x^\star\| + \|z_k - z^\star\| + \|\mathrm{prox}_{\gamma R}(0)\|\big) \to 0. \end{aligned}$$ □

We present below the proof of the main energy estimation in Section 3.

**Proof of Lemma 3.4.** The non-negativity of $\mathcal{D}_\Phi^{v^\star}(u_k)$ is rather obvious, owing to the convexity of $\Phi$. So, in the next, we focus on the second claim of the lemma. Define $y^{V^\perp} \stackrel{\text{def}}{=} \mathrm{P}_{V^\perp}(y), u_k^{V^\perp} \stackrel{\text{def}}{=} \mathrm{P}_{V^\perp}(u_k), z_k^{V^\perp} \stackrel{\text{def}}{=} \mathrm{P}_{V^\perp}(z_k)$ the projections of $y, u_k, z_k$ onto $V^\perp$ respectively.

From the update of $u_k$ in (1.6) and the definition of proximity operator, we get

$$(2x_k - z_k - u_{k+1})/\gamma_k - \nabla G(x_k) \in \partial R(u_{k+1}).$$

For the convexity of $R$, we obtain that, for every $y \in \mathcal{H}$,

$$\begin{aligned} R(y) &\geq R(u_{k+1}) + \langle(2x_k - z_k - u_{k+1})/\gamma_k - \nabla G(x_k),\, y - u_{k+1}\rangle \\ &= R(u_{k+1}) + \frac{1}{\gamma_k}\langle 2x_k - z_k - u_{k+1},\, y - u_{k+1}\rangle - \langle\nabla G(x_k),\, y - u_{k+1}\rangle. \end{aligned} \tag{A.7}$$

Notice that $u_{k+1} = x_k + z_{k+1} - z_k$. Then, the first inner product of the last line of (A.7) can be re-written as

$$\begin{aligned} &\langle 2x_k - z_k - u_{k+1},\, y - u_{k+1}\rangle \\ &= \langle x_k - z_k + z_k - z_{k+1},\, y - x_k - (z_{k+1} - z_k)\rangle \\ &= \langle x_k - z_k,\, y - x_k\rangle + \langle x_k - z_k,\, z_k - z_{k+1}\rangle - \langle x_k - y,\, z_k - z_{k+1}\rangle + \|z_{k+1} - z_k\|^2 \\ &= \langle x_k - z_k,\, y - x_k\rangle + \langle y - z_k,\, z_k - z_{k+1}\rangle + \|z_{k+1} - z_k\|^2 \\ &= \langle \mathrm{P}_V(z_k) - z_k,\, y - x_k\rangle + \tfrac{1}{2}\|z_{k+1} - y\|^2 - \tfrac{1}{2}\|z_{k+1} - z_k\|^2 - \tfrac{1}{2}\|z_k - y\|^2 + \|z_{k+1} - z_k\|^2 \\ &= -\langle z_k^{V^\perp},\, y\rangle + \tfrac{1}{2}\big(\|z_{k+1} - z_k\|^2 + \|z_{k+1} - y\|^2 - \|z_k - y\|^2\big), \end{aligned} \tag{A.8}$$

In the latter, we applied to $\langle y - z_k,\, z_k - z_{k+1}\rangle$ the usual Pythagoras relation, namely

$$2\langle c_2 - c_1,\, c_1 - c_3\rangle = \|c_2 - c_3\|^2 - \|c_1 - c_2\|^2 - \|c_1 - c_3\|^2.$$

Combining (A.8) with (A.7), we obtain

$$R(u_{k+1}) - R(y) \leq \langle\nabla G(x_k),\, y - u_{k+1}\rangle + \frac{1}{\gamma_k}\langle z_k^{V^\perp},\, y\rangle + \frac{1}{2\gamma_k}\big(\|z_k - y\|^2 - \|z_{k+1} - y\|^2 - \|z_{k+1} - z_k\|^2\big). \tag{A.9}$$

Since $G$ is convex, given any $x_k$ and $y \in \mathcal{H}$, we have

$$G(x_k) - G(y) \le \langle \nabla G(x_k),\, x_k - y\rangle. \tag{A.10}$$

Summing up (A.9) and (A.10) and rearranging the terms, we get

$$\begin{aligned}&\big(R(u_{k+1}) + G(x_k)\big) - \big(R(y) + G(y)\big) + \frac{1}{2\gamma_k}\big(\|z_{k+1} - y\|^2 - \|z_k - y\|^2\big) - \frac{1}{\gamma_k}\langle z_k^{V^\perp},\, y\rangle \\ &\le -\frac{1}{2\gamma_k}\|z_{k+1} - z_k\|^2 + \langle \nabla G(x_k),\, x_k - u_{k+1}\rangle.\end{aligned}$$

Since $G$ has Lipschitz continuous gradient, applying Lemma 2.1 yields

$$G(u_{k+1}) - G(x_k) \le \langle \nabla G(x_k),\, u_{k+1} - x_k\rangle + \frac{1}{2\beta_V}\|u_{k+1} - x_k\|^2.$$

Now sum up the above two inequalities and recall that $\xi_{k+1} \stackrel{\text{def}}{=} \frac{|\gamma - \beta_V|}{2\underline{\gamma}\beta_V}\|z_{k+1} - z_k\|^2$, to obtain

$$\begin{aligned}&\big(R(u_{k+1}) + G(u_{k+1})\big) - \big(R(y) + G(y)\big) + \frac{1}{2\gamma_k}\big(\|z_{k+1} - y\|^2 - \|z_k - y\|^2\big) - \frac{1}{\gamma_k}\langle z_k^{V^\perp},\, y\rangle \\ &\le -\frac{1}{2\gamma_k}\|z_{k+1} - z_k\|^2 + \frac{1}{2\beta_V}\|u_{k+1} - x_k\|^2 = \frac{\gamma_k - \beta_V}{2\gamma_k\beta_V}\|z_{k+1} - z_k\|^2 \le \frac{|\gamma_k - \beta_V|}{2\gamma_k\beta_V}\|z_{k+1} - z_k\|^2 \le \xi_{k+1}.\end{aligned} \tag{A.11}$$

Note that we applied again the equivalence $u_{k+1} = x_k + z_{k+1} - z_k$. Furthermore, define $\zeta_{k+1}^y \stackrel{\text{def}}{=} \frac{|\gamma_{k+1} - \gamma_k|}{2\underline{\gamma}^2}\|z_{k+1} - y\|^2$. Then, from (A.11), we have

$$\begin{aligned}\Phi(u_{k+1}) + \frac{1}{2\gamma_{k+1}}\|z_{k+1} - y\|^2 &= \Phi(u_{k+1}) + \frac{1}{2\gamma_k}\|z_{k+1} - y\|^2 + \big(\frac{1}{2\gamma_{k+1}} - \frac{1}{2\gamma_k}\big)\|z_{k+1} - y\|^2 \\ &\le \Phi(y) + \frac{1}{\gamma_k}\langle z_k^{V^\perp},\, y\rangle + \frac{1}{2\gamma_k}\|z_k - y\|^2 + \xi_{k+1} + \zeta_{k+1}^y \\ &= \Phi(y) + \frac{1}{\gamma_k}\langle z_k^{V^\perp},\, y^{V^\perp}\rangle + \frac{1}{2\gamma_k}\|z_k - y\|^2 + \xi_{k+1} + \zeta_{k+1}^y.\end{aligned} \tag{A.12}$$

Recall that $x_k \in V$. Hence, $\mathrm{P}_{V^\perp}(x_k) = 0$. Then, using (A.12), we have the following estimate for the Bregman divergence (defined in (3.4)):

$$\begin{aligned}&\mathcal{D}_\Phi^{v^\star}(u_{k+1}) - \mathcal{D}_\Phi^{v^\star}(y) = \Phi(u_{k+1}) - \Phi(x^\star) - \langle v^\star,\, u_{k+1}^{V^\perp} - y^{V^\perp}\rangle \\ &\le \frac{1}{\gamma_k}\langle z_k^{V^\perp},\, y^{V^\perp}\rangle - \langle v^\star,\, u_{k+1}^{V^\perp} - y^{V^\perp}\rangle + \frac{1}{2\gamma_k}\|z_k - y\|^2 - \frac{1}{2\gamma_{k+1}}\|z_{k+1} - y\|^2 + \xi_{k+1} + \zeta_{k+1}^y \\ &= \frac{1}{\gamma_k}\langle z_k^{V^\perp},\, y^{V^\perp}\rangle - \langle v^\star,\, x_k^{V^\perp} - y^{V^\perp}\rangle - \langle v^\star,\, z_{k+1}^{V^\perp}\rangle + \langle v^\star,\, z_k^{V^\perp}\rangle \\ &\qquad + \frac{1}{2\gamma_k}\|z_k - y\|^2 - \frac{1}{2\gamma_{k+1}}\|z_{k+1} - y\|^2 + \xi_{k+1} + \zeta_{k+1}^y \\ &= \frac{1}{2\gamma_k}\big(2\langle z_k^{V^\perp},\, y^{V^\perp}\rangle + 2\gamma_k\langle v^\star,\, y^{V^\perp}\rangle + 2\gamma_k\langle v^\star,\, z_k^{V^\perp}\rangle + \|z_k - y\|^2\big) \\ &\qquad + \frac{1}{2\gamma_{k+1}}\big(-2\gamma_{k+1}\langle v^\star,\, z_{k+1}^{V^\perp}\rangle - \|z_{k+1} - y\|^2\big) + \xi_{k+1} + \zeta_{k+1}^y \\ &= \frac{1}{2\gamma_k}\Big(-\big(\|z_k^{V^\perp} - y^{V^\perp}\|^2 - \|z_k^{V^\perp}\|^2 - \|y^{V^\perp}\|^2\big) + \big(\|y^{V^\perp} + \gamma_k v^\star\|^2 - \|y^{V^\perp}\|^2 - \|\gamma_k v^\star\|^2\big) \\ &\qquad + \big(\|z_k^{V^\perp} + \gamma_k v^\star\|^2 - \|z_k^{V^\perp}\|^2 - \|\gamma_k v^\star\|^2\big) + \|z_k - y\|^2\Big) \\ &\qquad + \frac{1}{2\gamma_{k+1}}\big(-\big(\|z_{k+1}^{V^\perp} + \gamma_{k+1}v^\star\|^2 - \|z_{k+1}^{V^\perp}\|^2 - \|\gamma_{k+1}v^\star\|^2\big) - \|z_{k+1} - y\|^2\big) + \xi_{k+1} + \zeta_{k+1}^y \\ &= \frac{1}{2\gamma_k}\big(\|y^{V^\perp} + \gamma_k v^\star\|^2 - 2\|\gamma_k v^\star\|^2 + \|z_k^{V^\perp} + \gamma_k v^\star\|^2 + \|z_k^V - y^V\|^2\big) \\ &\qquad + \frac{1}{2\gamma_{k+1}}\big(-\|z_{k+1}^{V^\perp} + \gamma_{k+1}v^\star\|^2 + \|z_{k+1}^{V^\perp}\|^2 + \|\gamma_{k+1}v^\star\|^2 - \|z_{k+1}^{V^\perp} - y^{V^\perp}\|^2 - \|z_{k+1}^V - y^V\|^2\big) + \xi_{k+1} + \zeta_{k+1}^y,\end{aligned}$$

where $y^V \overset{\text{def}}{=} \mathrm{P}_V(y), z_k^V \overset{\text{def}}{=} \mathrm{P}_V(z_k)$ are the projections of $y, z_k$ onto $V$ respectively. From above one, we deduce the following result

$$
\begin{aligned}
&\mathcal{D}_\Phi^{v^\star}(u_{k+1}) - \mathcal{D}_\Phi^{v^\star}(y) + \phi_{k+1} - \phi_k \\
&\leq \frac{1}{2\gamma_k}\big(\|y^{V^\perp} - \gamma_k v^\star\|^2 - 2\|\gamma_k v^\star\|^2\big) + \frac{1}{2\gamma_{k+1}}\big(\|z_{k+1}^{V^\perp}\|^2 + \|\gamma_{k+1}v^\star\|^2 - \|z_{k+1}^{V^\perp} - y^{V^\perp}\|^2\big) + \xi_{k+1} + \zeta_{k+1}^y \\
&= \frac{1}{2\gamma_k}\big(\|y^{V^\perp}\|^2 - 2\gamma_k\langle y^{V^\perp},\, v^\star\rangle - \|\gamma_k v^\star\|^2\big) + \frac{1}{2\gamma_{k+1}}\big(\|\gamma_{k+1}v^\star\|^2 - \|y^{V^\perp}\|^2 + 2\langle z_{k+1}^{V^\perp},\, y^{V^\perp}\rangle\big) + \xi_{k+1} + \zeta_{k+1}^y \\
&= \frac{1}{2\gamma_k}\|y^{V^\perp}\|^2 - \langle y^{V^\perp},\, v^\star\rangle - \frac{\gamma_k}{2}\|v^\star\|^2 + \frac{\gamma_{k+1}}{2}\|v^\star\|^2 - \frac{1}{2\gamma_{k+1}}\|y^{V^\perp}\|^2 + \frac{1}{\gamma_{k+1}}\langle z_{k+1}^{V^\perp},\, y^{V^\perp}\rangle + \xi_{k+1} + \zeta_{k+1}^y \\
&= \frac{\gamma_{k+1} - \gamma_k}{2\gamma_k\gamma_{k+1}}\|y^{V^\perp}\|^2 + \frac{\gamma_{k+1} - \gamma_k}{2}\|v^\star\|^2 + \frac{1}{\gamma_{k+1}}\langle z_{k+1}^{V^\perp} - \gamma_{k+1}v^\star,\, y^{V^\perp}\rangle + \xi_{k+1} + \zeta_{k+1}^y.
\end{aligned}
$$

In particular, taking $y = x^\star \in V$ in the last inequality and using the fact that $\mathrm{P}_{V^\perp}(x^\star) = 0$ lead to the desired result. □

The following lemma is very classical, see e.g. [24, Theorem 3.3.1].

**Lemma A.1.** *Let sequence $\{a_k\}_{k\in\mathbb{N}}$ be nonnegative, non-increasing and summable. Then $a_k = o(k^{-1})$.*

**Proof of Theorem 3.6.** Define $\theta_k \overset{\text{def}}{=} \min_{0\leq i\leq k} \mathcal{D}_\Phi^{v^\star}(u_i) \leq \mathcal{D}_\Phi^{v^\star}(u_k)$. Summing up the inequality (3.5) up to some $k \in \mathbb{N}$ yields

$$(k+1)\theta_k \leq \textstyle\sum_{i=0}^{k}\mathcal{D}_\Phi^{v^\star}(u_i) \leq \phi_0 + \frac{\gamma_\infty - \gamma_0}{2}\|v^\star\|^2 + \sum_{k\in\mathbb{N}}\xi_k + \sum_{k\in\mathbb{N}}\zeta_k.$$

Since $v^\star$ is bounded, so is $\phi_0$. Then, owing to Theorem 3.2, we have

$$\textstyle\sum_{k\in\mathbb{N}}\xi_k = \frac{|\gamma - \beta_V|}{2\underline{\gamma}\beta_V}\sum_{k\in\mathbb{N}}\|z_k - z_{k-1}\|^2 < +\infty.$$

Lastly, as $\{z_k\}_{k\in\mathbb{N}}$ is bounded, so is $\{\|z_k - x^\star\|\}_{k\in\mathbb{N}}$. Recall that, by assumptions, $\{\gamma_k\}_{k\in\mathbb{N}}$ converges to some $\gamma \in ]0, 2\beta_V[$ with $\{|\gamma_k - \gamma|\}_{k\in\mathbb{N}}$ being summable. Then we have that

$$
\begin{aligned}
\textstyle\sum_{k\in\mathbb{N}}\zeta_k &\leq \frac{1}{2\underline{\gamma}^2}\sup_{k\in\mathbb{N}}\|z_k - x^\star\|^2\textstyle\sum_{k\in\mathbb{N}}|\gamma_{k+1} - \gamma_k| \\
&\leq \frac{1}{2\underline{\gamma}^2}\sup_{k\in\mathbb{N}}\|z_k - x^\star\|^2\textstyle\sum_{k\in\mathbb{N}}\big(|\gamma_{k+1} - \gamma| + |\gamma_k - \gamma|\big) < +\infty.
\end{aligned}
$$

Summing up the above results, we have that $(k+1)\theta_k \leq C < +\infty$ holds for all $k \in \mathbb{N}$, which means

$$\theta_k = O(1/(k+1)).$$

Now, owing to the definition of $\theta_k$, we have

$$\textstyle\sum_{k\in\mathbb{N}}\theta_k \leq \sum_{k\in\mathbb{N}}\mathcal{D}_\Phi^{v^\star}(u_k) \leq \phi_0 + \frac{\gamma_\infty - \gamma_0}{2}\|v^\star\|^2 + \sum_{k\in\mathbb{N}}(\xi_k + \zeta_k) < +\infty.$$

Moreover, it is immediate that, for every $k \geq 1$,

$$\theta_k = \min(\mathcal{D}_\Phi^{v^\star}(u_k), \theta_{k-1}) \leq \theta_{k-1},$$

that is, the sequence $\{\theta_k\}_{k\in\mathbb{N}}$ is non-increasing. Invoking Lemma A.1 on $\{\theta_k\}_{k\in\mathbb{N}}$ concludes the proof.

For the ergodic rate, we start again from (3.5) and apply Jensen's inequality to $\mathcal{D}_\Phi^{v^\star}$ which is a convex function, and get

$$(k+1)\mathcal{D}_\Phi^{v^\star}(\bar{u}_k) \leq \textstyle\sum_{i=0}^{k}\mathcal{D}_\Phi^{v^\star}(u_i),$$

where the right-hand side is bounded by arguing as above. □

**Proof of Corollary 3.8.** First, weak convergence of the non-stationary FB iteration follows from Theorem 3.2.

On the one hand, specializing the inequality (3.5) to the case of FB splitting, we get

$$\Phi(x_{k+1})-\Phi(x^\star) \leq \frac{1}{2\gamma_k}\|x_k-x^\star\|^2-\frac{1}{2\gamma_{k+1}}\|x_{k+1}-x^\star\|^2+\frac{|\gamma-\beta|}{2\underline{\gamma}\beta}\|x_k-x_{k-1}\|^2+\frac{|\gamma_{k+1}-\gamma_k|}{2\underline{\gamma}^2}\|x_{k+1}-x^\star\|^2, \quad \text{(A.13)}$$

which means that

$$\begin{aligned}&\textstyle\sum_{k\in\mathbb{N}}\left(\Phi(x_k)-\Phi(x^\star)\right)\\ &\leq \frac{1}{2\gamma_0}\|x_0-x^\star\|^2+\frac{|\gamma-\beta|}{2\underline{\gamma}\beta}\textstyle\sum_{k\in\mathbb{N}}\|x_k-x_{k-1}\|^2+\frac{1}{\underline{\gamma}^2}\sup_{k\in\mathbb{N}}\|x_k-x^\star\|^2\textstyle\sum_{k\in\mathbb{N}}|\gamma_k-\gamma|<+\infty.\end{aligned}$$

On the other hand, by virtue of the inequality (A.11) in the proof of Lemma 3.4, given any $y\in\mathcal{H}$, we have

$$\Phi(x_{k+1})+\frac{1}{2\gamma_k}\|x_{k+1}-y\|^2 \leq \Phi(y)+\frac{1}{2\gamma_k}\|x_k-y\|^2+\big(\frac{1}{2\beta}-\frac{1}{2\gamma_k}\big)\|x_{k+1}-x_k\|^2.$$

Choosing $y=x_k$, we obtain

$$\begin{aligned}\big(\Phi(x_{k+1})-\Phi(x^\star)\big)-\big(\Phi(x_k)-\Phi(x^\star)\big) &\leq \big(\frac{1}{2\beta}-\frac{1}{\gamma_k}\big)\|x_{k+1}-x_k\|^2\\ &\leq \big(\frac{1}{2\beta}-\frac{1}{\gamma_k}\big)\|x_{k+1}-x_k\|^2\\ &\leq \big(\frac{1}{2\beta}-\frac{1}{\overline{\gamma}}\big)\|x_{k+1}-x_k\|^2=-\delta\|x_{k+1}-x_k\|^2,\end{aligned}$$

where $\delta=\frac{1}{\overline{\gamma}}-\frac{1}{2\beta}>0$ since $\overline{\gamma}<2\beta$. This implies that the sequence $\{\Phi(x_k)-\Phi(x^\star)\}_{k\in\mathbb{N}}$ is positive and non-increasing.

Summing up both sides of the above inequality and applying Lemma A.1 leads to the claimed result. □

# B Proofs of Section 4

## B.1 Riemannian Geometry

Let $\mathcal{M}$ be a $C^2$-smooth embedded submanifold of $\mathbb{R}^n$ around a point $x$. With some abuse of terminology, we shall state $C^2$-manifold instead of $C^2$-smooth embedded submanifold of $\mathbb{R}^n$. The natural embedding of a submanifold $\mathcal{M}$ into $\mathbb{R}^n$ permits to define a Riemannian structure and to introduce geodesics on $\mathcal{M}$, and we simply say $\mathcal{M}$ is a Riemannian manifold. We denote respectively $\mathcal{T}_{\mathcal{M}}(x)$ and $\mathcal{N}_{\mathcal{M}}(x)$ the tangent and normal space of $\mathcal{M}$ at point near $x$ in $\mathcal{M}$.

**Exponential map** Geodesics generalize the concept of straight lines in $\mathbb{R}^n$, preserving the zero acceleration characteristic, to manifolds. Roughly speaking, a geodesic is locally the shortest path between two points on $\mathcal{M}$. We denote by $\mathfrak{g}(t;x,h)$ the value at $t\in\mathbb{R}$ of the geodesic starting at $\mathfrak{g}(0;x,h)=x\in\mathcal{M}$ with velocity $\dot{\mathfrak{g}}(t;x,h)=\frac{d\mathfrak{g}}{dt}(t;x,h)=h\in\mathcal{T}_{\mathcal{M}}(x)$ (which is uniquely defined). For every $h\in\mathcal{T}_{\mathcal{M}}(x)$, there exists an interval $I$ around $0$ and a unique geodesic $\mathfrak{g}(t;x,h):I\to\mathcal{M}$ such that $\mathfrak{g}(0;x,h)=x$ and $\dot{\mathfrak{g}}(0;x,h)=h$. The mapping

$$\mathrm{Exp}_x:\mathcal{T}_{\mathcal{M}}(x)\to\mathcal{M},\ \ h\mapsto\mathrm{Exp}_x(h)=\mathfrak{g}(1;x,h),$$

is called *Exponential map*. Given $x,x'\in\mathcal{M}$, the direction $h\in\mathcal{T}_{\mathcal{M}}(x)$ we are interested in is such that

$$\mathrm{Exp}_x(h)=x'=\mathfrak{g}(1;x,h).$$

**Parallel translation** Given two points $x,x'\in\mathcal{M}$, let $\mathcal{T}_{\mathcal{M}}(x),\mathcal{T}_{\mathcal{M}}(x')$ be their corresponding tangent spaces. Define

$$\tau:\mathcal{T}_{\mathcal{M}}(x)\to\mathcal{T}_{\mathcal{M}}(x'),$$

the parallel translation along the unique geodesic joining $x$ to $x'$, which is isomorphism and isometry w.r.t. the Riemannian metric.

**Riemannian gradient and Hessian** For a vector $v \in \mathcal{N}_{\mathcal{M}}(x)$, the Weingarten map of $\mathcal{M}$ at $x$ is the operator $\mathfrak{W}_x(\cdot, v) : \mathcal{T}_{\mathcal{M}}(x) \to \mathcal{T}_{\mathcal{M}}(x)$ defined by

$$\mathfrak{W}_x(\cdot, v) = -\mathrm{P}_{\mathcal{T}_{\mathcal{M}}(x)} \mathrm{d}V[h],$$

where $V$ is any local extension of $v$ to a normal vector field on $\mathcal{M}$. The definition is independent of the choice of the extension $V$, and $\mathfrak{W}_x(\cdot, v)$ is a symmetric linear operator which is closely tied to the second fundamental form of $\mathcal{M}$, see [12, Proposition II.2.1].

Let $J$ be a real-valued function which is $C^2$ along the $\mathcal{M}$ around $x$. The covariant gradient of $J$ at $x' \in \mathcal{M}$ is the vector $\nabla_{\mathcal{M}} J(x') \in \mathcal{T}_{\mathcal{M}}(x')$ defined by

$$\langle \nabla_{\mathcal{M}} J(x'),\, h \rangle = \frac{d}{dt} J\big(\mathrm{P}_{\mathcal{M}}(x' + th)\big)\big|_{t=0}, \ \forall h \in \mathcal{T}_{\mathcal{M}}(x'),$$

where $\mathrm{P}_{\mathcal{M}}$ is the projection operator onto $\mathcal{M}$. The covariant Hessian of $J$ at $x'$ is the symmetric linear mapping $\nabla^2_{\mathcal{M}} J(x')$ from $\mathcal{T}_{\mathcal{M}}(x')$ to itself which is defined as

$$\langle \nabla^2_{\mathcal{M}} J(x')h,\, h \rangle = \frac{d^2}{dt^2} J\big(\mathrm{P}_{\mathcal{M}}(x' + th)\big)\big|_{t=0}, \ \forall h \in \mathcal{T}_{\mathcal{M}}(x'). \tag{B.1}$$

This definition agrees with the usual definition using geodesics or connections [36]. Now assume that $\mathcal{M}$ is a Riemannian embedded submanifold of $\mathbb{R}^n$, and that a function $J$ has a $C^2$-smooth restriction on $\mathcal{M}$. This can be characterized by the existence of a $C^2$-smooth extension (representative) of $J$, *i.e.* a $C^2$-smooth function $\widetilde{J}$ on $\mathbb{R}^n$ such that $\widetilde{J}$ agrees with $J$ on $\mathcal{M}$. Thus, the Riemannian gradient $\nabla_{\mathcal{M}} J(x')$ is also given by

$$\nabla_{\mathcal{M}} J(x') = \mathrm{P}_{\mathcal{T}_{\mathcal{M}}(x')} \nabla \widetilde{J}(x'), \tag{B.2}$$

and $\forall h \in \mathcal{T}_{\mathcal{M}}(x')$, the Riemannian Hessian reads

$$\begin{aligned} \nabla^2_{\mathcal{M}} J(x')h &= \mathrm{P}_{\mathcal{T}_{\mathcal{M}}(x')} \mathrm{d}(\nabla_{\mathcal{M}} J)(x')[h] = \mathrm{P}_{\mathcal{T}_{\mathcal{M}}(x')} \mathrm{d}\big(x' \mapsto \mathrm{P}_{\mathcal{T}_{\mathcal{M}}(x')} \nabla_{\mathcal{M}} \widetilde{J}\big)[h] \\ &= \mathrm{P}_{\mathcal{T}_{\mathcal{M}}(x')} \nabla^2 \widetilde{J}(x')h + \mathfrak{W}_{x'}\big(h, \mathrm{P}_{\mathcal{N}_{\mathcal{M}}(x')} \nabla \widetilde{J}(x')\big), \end{aligned} \tag{B.3}$$

where the last equality comes from [1, Theorem 1]. When $\mathcal{M}$ is an affine or linear subspace of $\mathbb{R}^n$, then obviously $\mathcal{M} = x + \mathcal{T}_{\mathcal{M}}(x)$, and $\mathfrak{W}_{x'}(h, \mathrm{P}_{\mathcal{N}_{\mathcal{M}}(x')} \nabla \widetilde{J}(x')) = 0$, hence (B.3) reduces to

$$\nabla^2_{\mathcal{M}} J(x') = \mathrm{P}_{\mathcal{T}_{\mathcal{M}}(x')} \nabla^2 \widetilde{J}(x') \mathrm{P}_{\mathcal{T}_{\mathcal{M}}(x')}.$$

See [25, 12] for more materials on differential and Riemannian manifolds.

We have the following proposition characterizing the parallel translation and the Riemannian Hessian of two close points in $\mathcal{M}$.

**Lemma B.1.** *Let $x, x'$ be two close points in $\mathcal{M}$, denote $\mathcal{T}_{\mathcal{M}}(x), \mathcal{T}_{\mathcal{M}}(x')$ be the tangent spaces of $\mathcal{M}$ at $x, x'$ respectively, and $\tau : \mathcal{T}_{\mathcal{M}}(x') \to \mathcal{T}_{\mathcal{M}}(x)$ be the parallel translation along the unique geodesic joining from $x$ to $x'$, then for the parallel translation we have, given any bounded vector $v \in \mathbb{R}^n$*

$$(\tau \mathrm{P}_{\mathcal{T}_{\mathcal{M}}(x')} - \mathrm{P}_{\mathcal{T}_{\mathcal{M}}(x)})v = o(\|v\|). \tag{B.4}$$

*The Riemannian Taylor expansion of $J \in C^2(\mathcal{M})$ at $x$ for $x'$ reads,*

$$\tau \nabla_{\mathcal{M}} J(x') = \nabla_{\mathcal{M}} J(x) + \nabla^2_{\mathcal{M}} J(x) \mathrm{P}_{\mathcal{T}_{\mathcal{M}}(x)}(x' - x) + o(\|x' - x\|). \tag{B.5}$$

**Proof.** See [30, Lemma B.1 and B.2]. □

**Lemma B.2.** *Let $\mathcal{M}$ be a $C^2$-smooth manifold, $\bar{x} \in \mathcal{M}$, $R \in \mathrm{PSF}_{\bar{x}}(\mathcal{M})$ and $\bar{u} \in \partial R(\bar{x})$. Let $\widetilde{R}$ be a smooth representative of $R$ on $\mathcal{M}$ near $x$, then given any $h \in T_{\bar{x}}$,*

(i) *when $\mathcal{M}$ is a general smooth manifold, if there holds $\bar{u} \in \mathrm{ri}\big(\partial R(\bar{x})\big)$, define the function $\overline{R}(x) = R(x) - \langle x,\, \bar{u} \rangle$, then*

$$\langle \mathrm{P}_{T_{\bar{x}}} \nabla^2_{\mathcal{M}} \overline{R}(\bar{x}) \mathrm{P}_{T_{\bar{x}}} h,\, h \rangle \geq 0. \tag{B.6}$$

(ii) *if $\mathcal{M}$ is affine/linear, then we have directly,*

$$\langle \mathrm{P}_{T_{\bar{x}}} \nabla^2_{\mathcal{M}} \widetilde{R}(\bar{x}) \mathrm{P}_{T_{\bar{x}}} h,\, h \rangle \geq 0. \tag{B.7}$$

**Proof.** See [30, Lemma 4.3]. □

## B.2 Proofs of main theorems

**Proof of Theorem 4.1.** From the updating of $u_{k+1}$ and the definition of proximity operator, we have that

$$(2x_k - z_k - u_{k+1})/\gamma_k - \nabla G(x_k) \in \partial R(u_{k+1}).$$

At convergence, we have

$$(x^\star - z^\star)/\gamma - \nabla G(x^\star) \in \partial R(x^\star).$$

Therefore,

$$\begin{aligned}
&\mathrm{dist}\big((x^\star - z^\star)/\gamma - \nabla G(x^\star), \partial R(u_{k+1})\big)\\
&\leq \|(2x_k - z_k - u_{k+1})/\gamma_k - (x^\star - z^\star)/\gamma - (\nabla G(x_k) - \nabla G(x^\star))\|\\
&\leq \|(2x_k - z_k - u_{k+1})/\gamma_k - (x^\star - z^\star)/\gamma_k + (x^\star - z^\star)/\gamma_k - (x^\star - z^\star)/\gamma\| + \|\nabla G(x_k) - \nabla G(x^\star)\|\\
&\leq \|(2x_k - z_k - u_{k+1})/\gamma_k - (x^\star - z^\star)/\gamma_k\| + \|(x^\star - z^\star)/\gamma_k - (x^\star - z^\star)/\gamma\| + \frac{1}{\beta_V}\|x_k - x^\star\|\\
&\leq \frac{1}{\gamma_k}\big(2\|x_k - x^\star\| + \|u_{k+1} - x^\star\| + \|z_k - z^\star\|\big) + \frac{|\gamma_k - \gamma|}{\gamma_k \gamma}\|(\mathrm{Id} - \mathrm{P}_V)(z^\star)\| + \frac{1}{\beta_V}\|x_k - x^\star\|\\
&\leq \frac{1}{\underline{\gamma}}\big(2\|x_k - x^\star\| + \|u_{k+1} - x^\star\| + \|z_k - z^\star\|\big) + \frac{|\gamma_k - \gamma|}{\underline{\gamma}^2}\|\mathrm{P}_{V^\perp}(z^\star)\| + \frac{1}{\beta_V}\|x_k - x^\star\|.
\end{aligned}$$

Theorem 3.2 allows to infer that the right hand side of the inequality converges to 0. In addition, since $R \in \Gamma_0(\mathbb{R}^n)$, $R$ is sub-differentially continuous at every point in its domain [42, Example 13.30], and in particular at $x^\star$. It then follows that $R(u_k) \to R(x^\star)$. Altogether, this shows that the conditions of [21, Theorem 5.3] are fulfilled for $R$, and the finite identification claim follows.

(a) In this case, $\mathcal{M}^R_{x^\star}$ is an affine subspace, *i.e.* $\mathcal{M}^R_{x^\star} = x^\star + T^R_{x^\star}$. Since $R$ is partly smooth at $x^\star$ relative to $\mathcal{M}^R_{x^\star}$, the sharpness property holds at all nearby points in $\mathcal{M}^R_{x^\star}$ [26, Proposition 2.10]. Thus for $k$ large enough, *i.e.* $u_k$ sufficiently close to $x^\star$ on $\mathcal{M}^R_{x^\star}$, we have indeed $\mathcal{T}_{u_k}(\mathcal{M}^R_{x^\star}) = T^R_{x^\star} = T^R_{u_k}$ as claimed.

(b) It is immediate to verify that a locally polyhedral function around $x^\star$ is indeed partly smooth relative to the affine subspace $x^\star + T^R_{x^\star}$, and thus, the first claim follows from (ii)(a). For the rest, it is sufficient to observe that by polyhedrality, for any $x \in \mathcal{M}^R_{x^\star}$ near $x^\star$, $\partial R(x) = \partial R(x^\star)$. Therefore, combining local normal sharpness [26, Proposition 2.10] and Lemma B.2 yields the second conclusion. □

**Proof of Theorem 4.4.** From (1.6), since $V$ is a subspace, then for $x_k$, we have

$$\begin{cases} x_k = \mathrm{P}_V(z_k), \\ x^\star = \mathrm{P}_V(z^\star), \end{cases} \iff \begin{cases} z_k - x_k \in \mathcal{N}_V(x_k), \\ z^\star - x^\star \in \mathcal{N}_V(x^\star). \end{cases}$$

Projecting onto $V$ leads to

$$x_k - x^\star = \mathrm{P}_V(z_k - z^\star). \tag{B.8}$$

Under the assumptions of Theorem 4.1, there exists $K \in \mathbb{N}$ large enough such that for all $k \geq K$, $u_k \in \mathcal{M}^R_{x^\star}$. Denote $T^R_{u_k}$ and $T^R_{x^\star}$ the tangent spaces corresponding to $u_k$ and $x^\star \in \mathcal{M}^R_{x^\star}$. Denote $\tau^R_k : T^R_{u_k} \to T^R_{x^\star}$ the parallel translation along the unique geodesic on $\mathcal{M}^R_{x^\star}$ joining $u_k$ to $x^\star$. Similar to (B.8), owing to [28, Lemma 5.1], we have for $u_k$ after identification that

$$u_k - x^\star = \mathrm{P}_{T^R_{x^\star}}(u_k - x^\star) + o(u_k - x^\star) \tag{B.9}$$

The update of $u_{k+1}$ in (1.6) and its convergence are respectively equivalent to

$$\begin{aligned}
2x_k - z_k - u_{k+1} - \gamma_k \nabla G(x_k) &\in \gamma_k \partial R(u_{k+1})\\
2x^\star - z^\star - x^\star - \gamma \nabla G(x^\star) &\in \gamma \partial R(x^\star).
\end{aligned}$$

Upon projecting onto the corresponding tangent spaces and applying the parallel translation $\tau_{k+1}$ from $u_{k+1}$ to $x^\star$, we get

$$\begin{aligned}\gamma_k\tau_{k+1}\nabla_{\mathcal{M}^R_{x^\star}} R(u_{k+1}) &= \tau_{k+1}\mathrm{P}_{T^R_{u_{k+1}}}\big(2x_k - z_k - u_{k+1} - \gamma_k\nabla G(x_k)\big)\\ &= \mathrm{P}_{T^R_{x^\star}}\big(2x_k - z_k - u_{k+1} - \gamma_k\nabla G(x_k)\big)\\ &\quad + \big(\tau_{k+1}\mathrm{P}_{T^R_{u_{k+1}}} - \mathrm{P}_{T^R_{x^\star}}\big)\big(2x_k - z_k - u_{k+1} - \gamma_k\nabla G(x_k)\big),\\ \gamma\nabla_{\mathcal{M}^R_{x^\star}} R(x^\star) &= \mathrm{P}_{T^R_{x^\star}}\big(2x^\star - z^\star - x^\star - \gamma\nabla G(x^\star)\big).\end{aligned}$$

Subtracting both equations, we obtain

$$\begin{aligned}&\gamma_k\tau_{k+1}\nabla_{\mathcal{M}^R_{x^\star}} R(u_{k+1}) - \gamma\nabla_{\mathcal{M}^R_{x^\star}} R(x^\star)\\ &= \gamma\tau_{k+1}\nabla_{\mathcal{M}^R_{x^\star}} R(u_{k+1}) - \gamma\nabla_{\mathcal{M}^R_{x^\star}} R(x^\star) + (\gamma_k-\gamma)\tau_{k+1}\nabla_{\mathcal{M}^R_{x^\star}} R(u_{k+1})\\ &= \mathrm{P}_{T^R_{x^\star}}\big((2x_k - z_k - u_{k+1} - \gamma_k\nabla G(x_k)) - (2x^\star - z^\star - x^\star - \gamma\nabla G(x^\star))\big)\\ &\quad + \underbrace{\big(\tau_{k+1}\mathrm{P}_{T^R_{u_{k+1}}} - \mathrm{P}_{T^R_{x^\star}}\big)(x^\star - z^\star - \gamma\nabla G(x^\star))}_{\textbf{Term 1}}\\ &\quad + \underbrace{\big(\tau_{k+1}\mathrm{P}_{T^R_{u_{k+1}}} - \mathrm{P}_{T^R_{x^\star}}\big)\big((2x_k - z_k - u_{k+1} - \gamma_k\nabla G(x_k)) - (2x^\star - z^\star - x^\star - \gamma\nabla G(x^\star))\big)}_{\textbf{Term 2}}.\end{aligned} \tag{B.10}$$

For the term $(\gamma_k-\gamma)\tau_{k+1}\nabla_{\mathcal{M}^R_{x^\star}} R(u_{k+1})$, since the Riemannian gradient $\nabla_{\mathcal{M}^R_{x^\star}} R(u_{k+1})$ is bounded on a bounded set, we have

$$(\gamma_k-\gamma)\tau_{k+1}\nabla_{\mathcal{M}^R_{x^\star}} R(u_{k+1}) = O(|\gamma_k-\gamma|). \tag{B.11}$$

For **Term 2**, owing to (B.4) and the boundedness of $\nabla G$, we have

$$\begin{aligned}&\big(\tau_{k+1}\mathrm{P}_{T^R_{u_{k+1}}} - \mathrm{P}_{T^R_{x^\star}}\big)\big((2x_k - z_k - u_{k+1} - \gamma_k\nabla G(x_k)) - (2x^\star - z^\star - x^\star - \gamma\nabla G(x^\star))\big)\\ &= o(\|(2x_k - z_k - u_{k+1} - \gamma_k\nabla G(x_k)) - (2x^\star - z^\star - x^\star - \gamma\nabla G(x^\star))\|)\\ &= o(\|(2x_k - z_k - u_{k+1} - \gamma\nabla G(x_k)) - (2x^\star - z^\star - x^\star - \gamma\nabla G(x^\star)) - (\gamma_k-\gamma)\nabla G(x_k)\|)\\ &= o(\|z_k - z^\star\|) + O(|\gamma_k-\gamma|).\end{aligned} \tag{B.12}$$

Then for **Term 1**, move to the other side of (B.10), combine the definition of $\overline{R}$ and the Riemannian Taylor expansion (B.5), we have

$$\begin{aligned}&\gamma\tau_{k+1}\nabla_{\mathcal{M}^R_{x^\star}} R(u_{k+1}) - \gamma\nabla_{\mathcal{M}^R_{x^\star}} R(x^\star) - \big(\tau_{k+1}\mathrm{P}_{T^R_{u_{k+1}}} - \mathrm{P}_{T^R_{x^\star}}\big)(x^\star - z^\star - \gamma\nabla G(x^\star))\\ &= \mathrm{P}_{T^R_{x^\star}}\nabla^2_{\mathcal{M}^R_{x^\star}}\overline{R}(x^\star)\mathrm{P}_{T^R_{x^\star}}(u_{k+1} - x^\star) + o(\|z_k - x^\star\|).\end{aligned} \tag{B.13}$$

Owing to [30, Lemma 4.3], we have that the Riemannian Hessian $\mathrm{P}_{T^R_{x^\star}}\nabla^2_{\mathcal{M}^R_{x^\star}}\overline{R}(x^\star)\mathrm{P}_{T^R_{x^\star}}$ is symmetric positive definite. For $\mathrm{P}_{T^R_{x^\star}}(\gamma_k\nabla G(x_k) - \gamma\nabla G(x^\star))$, since we assume that $F$ is locally $C^2$ around $x^\star$, then apply the Taylor expansion,

$$\begin{aligned}\gamma_k\nabla G(x_k) - \gamma\nabla G(x^\star) &= \gamma(\nabla G(x_k) - \nabla G(x^\star)) + (\gamma_k-\gamma)\nabla G(x_k)\\ &= \mathrm{P}_V\big(\nabla F(x_k) - \nabla F(x^\star)\big) + O(|\gamma_k-\gamma|)\\ &= \mathrm{P}_V\nabla^2 F(x_k - x^\star) + o(\|x_k - z^\star\|) + O(|\gamma_k-\gamma|)\\ &= \mathrm{P}_V\nabla^2 F\mathrm{P}_V(x_k - x^\star) + o(\|x_k - z^\star\|) + O(|\gamma_k-\gamma|)\\ &= \mathrm{P}_V\nabla^2 F\mathrm{P}_V(z_k - z^\star) + o(\|z_k - z^\star\|) + O(|\gamma_k-\gamma|).\end{aligned}$$

Then for (B.10) we have, recall that $H_{\overline{R}} \stackrel{\text{def}}{=} \mathrm{P}_{T^R_{x^\star}} \nabla^2_{\mathcal{M}^R_{x^\star}} \overline{R}(x^\star)\mathrm{P}_{T^R_{x^\star}}$ and $H_G \stackrel{\text{def}}{=} \mathrm{P}_V \nabla^2 F \mathrm{P}_V$,

$$\begin{aligned}
& H_{\overline{R}}(u_{k+1} - x^\star) = 2\mathrm{P}_{T^R_{x^\star}}(x_k - x^\star) - \mathrm{P}_{T^R_{x^\star}}(z_k - z^\star) - \mathrm{P}_{T^R_{x^\star}}(u_{k+1} - x^\star) \\
& \qquad - \gamma H_G(z_k - z^\star) + o(\|z_k - z^\star\|) + O(|\gamma_k - \gamma|) \\
\Longrightarrow \quad & (\mathrm{Id} + H_{\overline{R}})\mathrm{P}_{T^R_{x^\star}}(u_{k+1} - x^\star) = 2\mathrm{P}_{T^R_{x^\star}}(x_k - x^\star) - \mathrm{P}_{T^R_{x^\star}}(z_k - z^\star) - \gamma H_G(z_k - z^\star) \\
& \qquad + o(\|z_k - z^\star\|) + O(|\gamma_k - \gamma|) \\
\Longrightarrow \quad & \mathrm{P}_{T^R_{x^\star}}(u_{k+1} - x^\star) = 2M_{\overline{R}}\mathrm{P}_V(z_k - z^\star) - M_{\overline{R}}(z_k - z^\star) - \gamma M_{\overline{R}} H_G(z_k - z^\star) \\
& \qquad + o(\|z_k - z^\star\|) + O(|\gamma_k - \gamma|) \\
\Longrightarrow \quad & u_{k+1} - x^\star = 2M_{\overline{R}}\mathrm{P}_V(z_k - z^\star) - M_{\overline{R}}(z_k - z^\star) - \gamma M_{\overline{R}} H_G(z_k - z^\star) \\
& \qquad + o(\|z_k - z^\star\|) + O(|\gamma_k - \gamma|),
\end{aligned} \tag{B.14}$$

where we used several times the relation (B.9).

Summing up (B.8) and (B.14), we get

$$\begin{aligned}
& (z_k + u_{k+1} - x_k) - (z^\star + x^\star - x^\star) \\
& = (z_k - z^\star) + (u_{k+1} - x^\star) - (x_k - x^\star) \\
& = (\mathrm{Id} + 2M_{\overline{R}}\mathrm{P}_V - M_{\overline{R}} - \mathrm{P}_V - \gamma M_{\overline{R}} H_G)(z_k - z^\star) + o(\|z_k - z^\star\|) + O(|\gamma_k - \gamma|) \\
& = \mathcal{M}_\gamma(z_k - z^\star) + o(\|z_k - z^\star\|) + O(|\gamma_k - \gamma|).
\end{aligned}$$

Hence for the non-stationary FDR iteration, we have

$$\begin{aligned}
z_{k+1} - z^\star &= (1 - \lambda_k)(z_k - z^\star) + \lambda_k\big((z_k + u_{k+1} - x_k) - (z^\star + x^\star - x^\star)\big) \\
&= (1 - \lambda_k)(z_k - z^\star) + \lambda_k\big(\mathcal{M}_\gamma(z_k - z^\star) + o(\|z_k - z^\star\|) + O(|\gamma_k - \gamma|)\big) \\
&= (1 - \lambda_k)(z_k - z^\star) + \lambda_k\mathcal{M}_\gamma(z_k - z^\star) + o(\|z_k - z^\star\|) + \chi_k \\
&= \mathcal{M}_{\gamma,\lambda}(z_k - z^\star) - (\lambda_k - \lambda)(\mathrm{Id} - \mathcal{M}_\gamma)(z_k - z^\star) + o(\|z_k - z^\star\|) + \chi_k.
\end{aligned}$$

Since we have

$$\lim_{k \to +\infty} \frac{\|(\lambda_k - \lambda)(\mathrm{Id} - \mathcal{M}_\gamma)(z_k - z^\star)\|}{\|z_k - z^\star\|} \leq \lim_{k \to +\infty} \frac{|\lambda_k - \lambda|\|\mathrm{Id} - \mathcal{M}_\gamma\|\|z_k - z^\star\|}{\|z_k - z^\star\|} = \lim_{k \to +\infty} |\lambda_k - \lambda|\|\mathrm{Id} - \mathcal{M}_\gamma\| = 0,$$

which means that

$$z_{k+1} - z^\star = \mathcal{M}_{\gamma,\lambda}(z_k - z^\star) + \psi_k + \chi_k,$$

and the claimed result is obtained. ☐

**Proof of Lemma 4.5.** Since $W_{\overline{R}}$ is firmly non-expansive by Lemma 4.3, it follows from [4, Example 4.7] that $M_{\overline{R}}$ is firmly non-expansive, hence $\mathcal{R}_{M_{\overline{R}}} \stackrel{\text{def}}{=} 2M_{\overline{R}} - \mathrm{Id}$ is non-expansive. Similarly, as $\mathrm{P}_V$ is firmly non-expansive, $\mathcal{R}_V \stackrel{\text{def}}{=} 2\mathrm{P}_V - \mathrm{Id}$ is non-expansive. As a result, we have $\frac{1}{2}\big(\mathcal{R}_{M_{\overline{R}}}\mathcal{R}_V + \mathrm{Id}\big)$ is firmly non-expansive [4, Proposition 4.21(i)-(ii)]. Then for $\mathrm{Id} - \gamma H_G$, given $\gamma \in [0, 2\beta_V]$ it is $\frac{2\beta_V}{4\beta_V - \gamma}$-averaged non-expansive. Therefore, owing to Lemma 3.1, we have the averaged property of $\mathcal{M}_\gamma$ and $\mathcal{M}_{\gamma,\lambda}$. We deduce from [4, Proposition 5.15] that $\mathcal{M}_\gamma$ and $\mathcal{M}_{\gamma,\lambda}$ are convergent, *i.e.* the power $\mathcal{M}^k_{\gamma,\lambda}$ exists as $k$ approaches $+\infty$ which is denoted as $\mathcal{M}^\infty_\gamma$. Moreover, $\mathcal{M}^k_{\gamma,\lambda} - \mathcal{M}^\infty_\gamma = (\mathcal{M}_{\gamma,\lambda} - \mathcal{M}^\infty_\gamma)^k$, $\forall k \in \mathbb{N}$, and $\rho(\mathcal{M}_{\gamma,\lambda} - \mathcal{M}^\infty_\gamma) < 1$ by [5, Theorem 2.12]. The second claim of the lemma is classical using the spectral radius formula (4.7), see *e.g.* [5, Theorem 2.12(i)]. ☐

**Proof of Theorem 4.6.**

(i) Let $K \in \mathbb{N}$ sufficiently large such that the locally linearized iteration (4.6) holds, then we have for $k \geq K$

$$\begin{aligned}
z_{k+1} - z^\star &= \mathcal{M}_{\gamma,\lambda}(z_k - z^\star) + \psi_k + \chi_k \\
&= \mathcal{M}_{\gamma,\lambda}\big(\mathcal{M}_{\gamma,\lambda}(z_{k-1} - z^\star) + \psi_{k-1} + \chi_{k-1}\big) + \psi_k + \chi_k \\
&= \mathcal{M}^{k+1-K}_{\gamma,\lambda}(z_K - z^\star) + \textstyle\sum_{j=K}^k \mathcal{M}^{k-j}_{\gamma,\lambda}(\psi_j + \chi_j).
\end{aligned} \tag{B.15}$$

Since $z_k \to z^\star$ and $\mathcal{M}_{\gamma,\lambda}$ is convergent to $\mathcal{M}_\gamma^\infty$ by Lemma 4.5, taking the limit as $k \to +\infty$, we have for all finite $p \geq K$,

$$\lim_{k\to+\infty} \sum_{j=p}^{k} \mathcal{M}_{\gamma,\lambda}^{k-j}(\psi_j + \chi_j) = -\mathcal{M}_\gamma^\infty(z_p - z^\star). \tag{B.16}$$

Using (B.16) in (B.15), we get

$$\begin{aligned}
&z_{k+1} - z^\star \\
&= (\mathcal{M}_{\gamma,\lambda} - \mathcal{M}_\gamma^\infty)(z_k - z^\star) + \psi_k + \chi_k - \lim_{l\to+\infty} \sum_{j=k}^{l} \mathcal{M}_{\gamma,\lambda}^{l-j}(\psi_j + \chi_j) \\
&= (\mathcal{M}_{\gamma,\lambda} - \mathcal{M}_\gamma^\infty)(z_k - z^\star) + \psi_k + \chi_k - \lim_{l\to+\infty} \sum_{j=k+1}^{l} \mathcal{M}_{\gamma,\lambda}^{l-j}(\psi_j + \chi_j) - \mathcal{M}_\gamma^\infty(\psi_k + \chi_k) \\
&= (\mathcal{M}_{\gamma,\lambda} - \mathcal{M}_\gamma^\infty)(z_k - z^\star) + (\mathrm{Id} - \mathcal{M}_\gamma^\infty)(\psi_j + \chi_j) + \mathcal{M}_\gamma^\infty(z_{k+1} - z^\star).
\end{aligned}$$

It is also immediate to see from Lemma 4.5 that $\|\mathrm{Id} - \mathcal{M}_\gamma^\infty\| \leq 1$ and

$$(\mathcal{M}_{\gamma,\lambda} - \mathcal{M}_\gamma^\infty)(\mathrm{Id} - \mathcal{M}_\gamma^\infty) = \mathcal{M}_{\gamma,\lambda} - \mathcal{M}_\gamma^\infty.$$

Rearranging the terms gives the claimed equivalence.

(ii) Under polyhedrality and constant parameters, we have from Theorem 4.4 that $o(\|z_k - z^\star\|)$ and $O(\lambda_k|\gamma_k - \gamma|)$ vanish, and the result follows. □

**Proof of Theorem 4.7.**

(i) Let $K \in \mathbb{N}$ sufficiently large such that (4.8) holds. We then have from Corollary 4.6(i)

$$\begin{aligned}
&(\mathrm{Id} - \mathcal{M}_\gamma^\infty)(z_{k+1} - z^\star) \\
&= (\mathcal{M}_{\gamma,\lambda} - \mathcal{M}_\gamma^\infty)^{k+1-K}(\mathrm{Id} - \mathcal{M}_\gamma^\infty)(z_K - z^\star) + \sum_{j=K}^{k} (\mathcal{M}_{\gamma,\lambda} - \mathcal{M}_\gamma^\infty)^{k-j}\big((\mathrm{Id} - \mathcal{M}_\gamma^\infty)\psi_j + \chi_j\big).
\end{aligned}$$

Since $\rho(\mathcal{M}_{\gamma,\lambda} - \mathcal{M}_\gamma^\infty) < 1$ by Lemma 4.5, from the spectral radius formula, we know that for every $\rho \in ]\rho(\mathcal{M}_{\gamma,\lambda} - \mathcal{M}_\gamma^\infty), 1[$, there is a constant $C$ such that

$$\|(\mathcal{M}_{\gamma,\lambda} - \mathcal{M}_\gamma^\infty)^j\| \leq C\rho^j$$

holds for all integers $j$. We thus get

$$\begin{aligned}
&\|(\mathrm{Id} - \mathcal{M}_\gamma^\infty)(z_{k+1} - z^\star)\| \\
&\leq C\rho^{k+1-K}\|z_K - z^\star\| + C\sum_{j=K}^{k} \rho^{k-j}\chi_j + C\sum_{j=K}^{k} \rho^{k-j}\|(\mathrm{Id} - \mathcal{M}_{\gamma,\lambda})\psi_j\| \\
&= C\rho^{k+1-K}\big(\|z_K - z^\star\| + \rho^{K-1}\sum_{j=K}^{k} \tfrac{\chi_j}{\rho^j}\big) + C\sum_{j=K}^{k} \rho^{k-j}\|(\mathrm{Id} - \mathcal{M}_{\gamma,\lambda})\psi_j\|,
\end{aligned} \tag{B.17}$$

By assumption, $\chi_j = C'\eta^j$, for some constant $C' \geq 0$ and $\eta < \rho$, and we have

$$\rho^{K-1}\sum_{j=K}^{k} \frac{\chi_j}{\rho^j} \leq C'\rho^{K-1}\sum_{j=K}^{\infty} (\eta/\rho)^j = \frac{C'\eta^K}{\rho - \eta} < +\infty.$$

Setting $C'' = C(\|z_K - z^\star\| + \frac{C'\eta^K}{\rho-\eta}) < +\infty$, we obtain

$$\|(\mathrm{Id} - \mathcal{M}_\gamma^\infty)(z_{k+1} - z^\star)\| \leq C''\rho^{k+1-K} + C\sum_{j=K}^{k} \rho^{k-j}\|(\mathrm{Id} - \mathcal{M}_\gamma^\infty)\psi_j\|.$$

This, together with the fact that $\|(\mathrm{Id} - \mathcal{M}_\gamma^\infty)\psi_j\| = o(\|(\mathrm{Id} - \mathcal{M}_\gamma^\infty)(z_j - z^\star)\|)$ yields the claimed result.

(ii) From Corollary 4.6, we have

$$z_k - z^\star = (\mathcal{M}_{\gamma,\lambda} - \mathcal{M}_\gamma^\infty)^{k+1-K}(z_K - z^\star),$$

hence the result follows. □

# References

[1] P-A. Absil, R. Mahony, and J. Trumpf. An extrinsic look at the Riemannian Hessian. In *Geometric Science of Information*, pages 361–368. Springer, 2013.

[2] H. Attouch and J. Peypouquet. The rate of convergence of Nesterov's accelerated Forward–Backward method is actually faster than $1/k^2$. *SIAM Journal on Optimization*, 26(3):1824–1834, 2016.

[3] J. B. Baillon and G. Haddad. Quelques propriétés des opérateurs angle-bornés etn-cycliquement monotones. *Israel Journal of Mathematics*, 26(2):137–150, 1977.

[4] H. Bauschke and P. L. Combettes. *Convex Analysis and Monotone Operator Theory in Hilbert Spaces*. Springer, 2011.

[5] H. H. Bauschke, J.Y. Bello Cruz, T.A. Nghia, H. M. Phan, and X. Wang. Optimal rates of convergence of matrices with applications. *Numerical Algorithms*, 2016. in press (arxiv:1407.0671).

[6] A. Beck and M. Teboulle. A fast iterative shrinkage-thresholding algorithm for linear inverse problems. *SIAM Journal on Imaging Sciences*, 2(1):183–202, 2009.

[7] D. P. Bertsekas. *Nonlinear programming*. Athena scientific Belmont, 1999.

[8] J. Bolte, T. P. Nguyen, J. Peypouquet, and B. W. Suter. From error bounds to the complexity of first-order descent methods for convex functions. *Mathematical Programming*, 165(2):471–507, Oct 2017.

[9] H. Brézis. *Opérateurs Maximaux Monotones et Semi-Groupes de Contractions dans les Espaces de Hilbert*. North-Holland/Elsevier, New York, 1973.

[10] L. M. Briceño-Arias. Forward-douglas–rachford splitting and forward-partial inverse method for solving monotone inclusions. *Optimization*, 64(5):1239–1261, 2015.

[11] A. Chambolle and C. Dossal. On the convergence of the iterates of the "fast iterative shrinkage/thresholding algorithm". *Journal of Optimization Theory and Applications*, 166(3):968–982, 2015.

[12] I. Chavel. *Riemannian geometry: a modern introduction*, volume 98. Cambridge University Press, 2006.

[13] P. L. Combettes. Quasi-Fejérian analysis of some optimization algorithms. *Studies in Computational Mathematics*, 8:115–152, 2001.

[14] P. L. Combettes and V. R. Wajs. Signal recovery by proximal forward-backward splitting. *Multiscale Modeling & Simulation*, 4(4):1168–1200, 2005.

[15] A. Daniilidis, D. Drusvyatskiy, and A. S. Lewis. Orthogonal invariance and identifiability. *SIAM Journal on Matrix Analysis and Applications*, 35(2):580–598, 2014.

[16] D. Davis. Convergence rate analysis of the Forward–Douglas–Rachford splitting scheme. *SIAM Journal on Optimization*, 25(3):1760–1786, 2015.

[17] D. Davis and W. Yin. A three-operator splitting scheme and its optimization applications. *Set-Valued and Variational Analysis*, 2017.

[18] J. Douglas and H. H. Rachford. On the numerical solution of heat conduction problems in two and three space variables. *Transactions of the American mathematical Society*, 82(2):421–439, 1956.

[19] D. Drusvyatskiy and A.D. Lewis. Error bounds, quadratic growth, and linear convergence of proximal methods. Technical Report arXiv:1602.06661, 2016.

[20] B. R. Gaines and H. Zhou. Algorithms for fitting the constrained lasso. *arXiv preprint arXiv:1611.01511*, 2016.

[21] W. L. Hare and A. S. Lewis. Identifying active constraints via partial smoothness and prox-regularity. *Journal of Convex Analysis*, 11(2):251–266, 2004.

[22] T. He. Lasso and general l1-regularized regression under linear equality and inequality constraints. 2011.

[23] G. M. James, C. Paulson, and P. Rusmevichientong. The constrained lasso. *Manuscript*, 2012.

[24] K. Knopp. *Theory and application of infinite series*. Courier Corporation, 2013.

[25] J. M. Lee. *Smooth manifolds*. Springer, 2003.

[26] A. S. Lewis. Active sets, nonsmoothness, and sensitivity. *SIAM Journal on Optimization*, 13(3):702–725, 2003.

[27] G. Li and T. K. Pong. Calculus of the exponent of Kurdyka–Łojasiewicz inequality and its applications to linear convergence of first-order methods. *Foundations of Computational Mathematics*, Aug 2017.

[28] J. Liang, J. Fadili, and G. Peyré. Local linear convergence of Forward–Backward under partial smoothness. In *Advances in Neural Information Processing Systems*, pages 1970–1978, 2014.

[29] J. Liang, J. Fadili, and G. Peyré. Convergence rates with inexact non-expansive operators. *Mathematical Programming: Series A*, 159(1-2):403–434, 2016.

[30] J. Liang, J. Fadili, and G. Peyré. Activity identification and local linear convergence of Forward–Backward-type methods. *SIAM Journal on Optimization*, 27(1):408–437, 2017.

[31] J. Liang, J. Fadili, and G. Peyré. Local convergence properties of Douglas–Rachford and Alternating Direction Method of Multipliers. *Journal of Optimization Theory and Applications*, 172(3):874–913, 2017.

[32] J. Liang, J. Fadili, and G. Peyré. Local linear convergence analysis of Primal–Dual splitting methods. *arXiv preprint arXiv:1705.01926*, 2017.

[33] P. L. Lions and B. Mercier. Splitting algorithms for the sum of two nonlinear operators. *SIAM Journal on Numerical Analysis*, 16(6):964–979, 1979.

[34] Z.-Q. Luo and P. Tseng. On the linear convergence of descent methods for convex essentially smooth minimization. *SIAM Journal on Control and Optimization*, 30(2):408–425, 1992.

[35] Z.-Q. Luo and P. Tseng. Error bounds and convergence analysis of feasible descent methods: a general approach. *Annals of Operations Research*, 46(1):157–178, Mar 1993.

[36] S. A. Miller and J. Malick. Newton methods for nonsmooth convex minimization: connections among-Lagrangian, Riemannian Newton and SQP methods. *Mathematical programming*, 104(2-3):609–633, 2005.

[37] Y. Nesterov. A method for solving the convex programming problem with convergence rate $O(1/k^2)$. *Dokl. Akad. Nauk SSSR*, 269(3):543–547, 1983.

[38] Y. Nesterov. *Introductory lectures on convex optimization: A basic course*, volume 87. Springer, 2004.

[39] N. Ogura and I. Yamada. Non-strictly convex minimization over the fixed point set of an asymptotically shrinking nonexpansive mapping. *Numerical Functional Analysis and Optimization*, 23(1-2):113–137, 2002.

[40] H. Raguet, M. J. Fadili, and G. Peyré. Generalized forward-backward splitting. *SIAM Journal on Imaging Sciences*, 6(3):1199–1226, 2013.

[41] R. T. Rockafellar. *Convex analysis*, volume 28. Princeton university press, 1997.

[42] R. T. Rockafellar and R. Wets. *Variational analysis*, volume 317. Springer Verlag, 1998.

[43] S. Vaiter, C. Deledalle, J. M. Fadili, G. Peyré, and C. Dossal. The degrees of freedom of partly smooth regularizers. *Annals of the Institute of Mathematical Statistics*, 2015. to appear.

[44] S. J. Wright. Identifiable surfaces in constrained optimization. *SIAM Journal on Control and Optimization*, 31(4):1063–1079, 1993.

[45] Z. Zhou and A. M.-C. So. A unified approach to error bounds for structured convex optimization problems. *Mathematical Programming*, 165(2):689–728, Oct 2017.